\documentclass[12pt]{article}
\usepackage{latexsym}
\usepackage{amsmath}
\usepackage{amssymb}
\usepackage{amsthm}
\usepackage{amscd}
\usepackage[dvipdfmx]{graphicx}
\input cyracc.def
\newfont{\tencyr}{wncyr10}

\begin{document}

\begin{center}
{\large \bf  Arithmetic topology in Ihara theory II:}\\
{\large \bf  Milnor invariants, dilogarithmic Heisenberg coverings and triple power residue symbols}
 \end{center}

 \vspace{.05cm}
 
\begin{center}
Hikaru HIRANO and Masanori MORISHITA
\end{center}

\begin{center}
{\em Dedicated to Professor Yasutaka Ihara}
\end{center}


{\small
Abstract: We introduce mod $l$ Milnor invariants of a Galois element associated to Ihara's Galois representation on the pro-$l$ fundamental group of a punctured projective line ($l$ being a prime number), as arithmetic analogues of Milnor invariants of a pure braid. We then show that triple quadratic (resp. cubic) residue symbols of primes in the rational (resp. Eisenstein) number field are expressed by mod $2$ (resp. mod $3$)  triple Milnor invariants of Frobenius elements. For this, we introduce dilogarithmic  mod $l$ Heisenberg ramified covering ${\cal D}^{(l)}$ of $\mathbb{P}^1$, which may be regarded as a higher analog of the dilogarithmic function, for the gerbe associated to the mod $l$ Heisenberg group, and we study the monodromy transformations of certain functions on  ${\cal D}^{(l)}$ along the pro-$l$ longitudes of Frobenius elements for $l=2,3$. }

\vspace{.8cm}

\begin{center}
 {\bf  Introduction}
 \end{center}
 
 In [KMT], following the analogy between the Artin representation of a pure braid group and the Ihara representation of the absolute Galois group ${\rm Gal}_k := {\rm Gal}(\overline{k}/k)$ of a number field $k$ on the pro-$l$ fundamental group of a punctured projective line ([I1],[I2]), $l$-adic Milnor invariants $\overline{\mu}^{(l)}(g;I) := \mu^{(l)}(g;I) \; \mbox{mod} \; \Delta(g;I) $ of each Galois element $g \in {\rm Gal}_k$ were introduced as arithmetic analogues of Milnor invariants of a pure braid ([MK; Chapter 6, 4],[Kd; 1.2]), where $l$ is a prime number, $I$ is a multi-index representing punctured points and $\Delta(g; I)$ is a certain indeterminacy (cf. Subsection 1.3). They were shown to enjoy some properties similar to those of Milnor invariants of pure braids. In principle, the information on the Ihara representation is encoded in $l$-adic Milnor numbers $\mu^{(l)}(g;I)$ for all $g$ and $I$.  

On the other hand, based on the analogies between knots and primes ([Mo4]), we have mod $2$ Milnor invariants $\mu_2(J)$ of certain rational primes ([Mo1]$\sim$[Mo4]), as arithmetic analogues of Milnor invariants of a link ([Mi1], [Mi2], [T]), where a multi-index $J$ represents an ordered set of primes (cf. Subsection 3.2). For example, $(-1)^{\mu_2(12)}$ coincides with the Legendre symbol $(p_1/p_2)$. Assuming $\mu_2(ij) = 0$ ($1\leq i,j\leq 3$), $(-1)^{\mu_2(123)}$ is proved to equal the triple quadratic residue symbol $[p_1,p_2,p_3]$ introduced by R\'{e}dei ([R]), which  describes the decomposition of $p_3$ in a certain dihedral extension, determined by $p_1$ and $p_2$,  of degree $8$ over $\mathbb{Q}$. Recently, mod $3$ Milnor invariants $\mu_3(ij)$ and $\mu_3(123)$ were introduced for certain primes $\frak{p}_i = (\pi_i)$ ($1 \leq i \leq 3$) of the Eisenstein number field $\mathbb{Q}(\zeta_3)$, $\zeta_3 := \exp(2\pi\sqrt{-1}/3)$ ([AMM]). As in the mod $2$ case, $\zeta_3^{\mu_3(12)}$ coincides with the cubic residue symbol $(\pi_1/\pi_2)_3$. Assuming $\mu_3(ij) = 0$ for $1\leq i, j \leq 3$, $\zeta_3^{\mu_3(123)}$ is proved to equal the triple cubic residue symbol $[\frak{p}_1,\frak{p}_2,\frak{p}_3]_3$, which  describes the decomposition of $\frak{p}_3$ in a certain mod $3$ Heisenberg extension,  determined by $\frak{p}_1$ and $\frak{p}_2$,  of degree $27$ over $\mathbb{Q}(\zeta_3)$ ([ibid]). We note that a key ingredient to introduce these Milnor invariants of primes is the theory of pro-$l$ extensions of number fields with restricted ramification due to Koch et al. (cf. [Kc]).

Since Milnor invariants of a braid $b$ coincide with those of the link obtained by closing $b$, the analogy with topology suggests to ask if there would be any relation between mod $l$ Milnor invariants $\mu_l(g;I) := \mu^{(l)}(g;I) \; \mbox{mod}\; l$ of Galois elements and mod $l$ Milnor invariants $\mu_l(J)$ of primes. This question may be of arithmetic interest and importance, because such a relation would reveal  a connection between Ihara theory and the classical arithmetic of pro-$l$ extensions of number fields. In this paper, we study this question.  Let us describe our results in the following.

We firstly interpret the pro-$l$ longitudes of a Galois element ([KMT; $\S 3.2$) in terms of certain pro-$l$ paths and show that  the $l$-th power residue symbol can be given by a mod $l$ Milnor invariant $\mu_l(\sigma;I)$ of a Frobenius element $\sigma$ with $|I| = 2$. Our main result is that triple quadratic (resp. cubic) residue symbols can be expressed by mod $2$ (resp. mod $3$) triple Milnor invariants of Frobenius elements (Theorem 4.1.10, Theorem 4.2.14), and hence answers the above question for triple Milnor invariants.  For this, 
we introduce a certain mod $l$ Heisenberg ramified covering ${\cal D}^{(l)}$ of $\mathbb{P}^1$ , called the {\em dilogarithmic mod $l$ Heisenberg ramified covering}, which may be regarded as a higher analog of the dilogarithmic function, for the gerbe associated to the mod $l$ Heisenberg group. We then study the monodromy transformations of certain functions on ${\cal D}^{(l)}$  along the pro-$l$ longitudes of Frobenius elements for $l = 2, 3$, to obtain our main result. Our method is closely related with Wojtkowiak's work ([NW], [W1] $\sim$ [W5]).

Here are the contents of this paper. In Section 1, we introduce the pro-$l$ longitudes of a Galois element in terms of pro-$l$ paths, and then introduce mod $l$ Milnor invariants of Galois elements. In Section 2, we introduce certain mod $l$ Heisenberg ramified coverings of $\mathbb{P}^1$, and explain the analogies with the dilogarithmic function. In  Section 3, we recall mod $2$ (resp. mod $3$) Milnor invariants of primes of $\mathbb{Q}$ (resp. $\mathbb{Q}(\zeta_3)$). In Section 4, we interpret the R\'{e}dei symbol and the triple cubic residue symbol by mod $l$ Milnor invariants of Frobenius elements for $l=2$ and $3$, respectively, by computing the monodromy transformations of certain functions on the dilogarithmic mod $l$ Heisenberg coverings in Section 2 along the pro-$l$ longitudes, and deduce the relations between mod $l$ Milnor invariants of Galois elements in Ihara theory and mod $l$ Milnor invariants of primes for $l = 2,3$.\\
\\
{\em Acknowledgement.} We would like to thank Yasushi Mizusawa, Hiroaki Nakamura, Yuji Terashima, Hiroshi Tsunogai and Zdzis\l aw Wojtkowiak for useful communications. 
 Especially, we thank Terashima for discussions on the subsection 2.2. We would like to thank the referee for useful comments which improved  the earlier version. The second author is partly supported by JSPS KAKENHI Grant Number JP17H02837, Grant-in-Aid for Scientific Research (B).\\
\\
{\em Notation.} Throughout this paper, $l$ denotes a prime number.\\
For a number field $K$, ${\cal O}_K$ denotes the ring of integers of $K$.\\
For subgroups $A, B$ of a topological group $G$, $[A,B]$ stands for the closed subgroup of $G$ generated by commutators  $[a,b] := aba^{-1}b^{-1}$ for $a \in A, b \in B$. \\

 \begin{center}
 {\bf 1. Mod $l$ Milnor invariants of Galois elements in Ihara theory}
 \end{center}
 
In [KMT; $\S 3$], following the analogy with Milnor invariants of braids associated to the Artin representation ([MK; Chapter 6, 4]), we introduced $l$-adic Milnor invariants
of each Galois element, in a group theoretic manner, as the Magnus coefficients of the pro-$l$ longitudes of a Galois element associated to the Ihara representation. In this section, following [I3], [NW] and Wojtkowiak's series of papers [W1]$\sim$[W4], we interpret the pro-$l$ longitudes  in terms of pro-$l$ paths and then introduce mod  $l$ Milnor invariants of a Galois element.  We show that mod $l$ Milnor invariants for indices of length 2 are given by $l$-th power residue symbols.\\
\\
 {\bf 1.1. The Ihara representation.} Let $k$ be a fixed finite algebraic number field in the field $\mathbb{C}$ of complex numbers and $\overline{k}$ a fixed algebraic closure of $k$ in $\mathbb{C}$. Let $\mathbb{P}^1$ be the projective $t$-line over $k$. Let $a_1, \dots, a_r$ be distinct $r$ numbers in $k$ ($r \geq 2$), identified with $k$-rational points on $\mathbb{P}^1$, and let $A := \{a_0, a_1, \dots , a_r \}$ with $a_0 = \infty$. We let $X := \mathbb{P}^1 \setminus A = {\rm Spec}\, k[t, (t-a_j)^{-1} (1\leq j \leq r)]$ and $X_{\overline{k}} := X \otimes_k \overline{k}$. For each $j$ $0 \leq j \leq r$, let $v_j$ be a $k$-rational tangential base point on $X$ at $a_j$ ([N; I]), which may be regarded as a tangential base point on the complex manifold $X(\mathbb{C})$ at $a_j$, a tangent vector on $\mathbb{P}^1(\mathbb{C})$ at $a_i$. Following [W1; $\S 2$], the geometric generators of $\pi_1(X(\mathbb{C}); v_0)$ are defined as follows. Let $x_0$ be a small circle on $\mathbb{P}^1(\mathbb{C})$ around $a_0 = \infty$ starting from $v_0$ in the opposite clockwise way. Choose a point $v_0' \in X(\mathbb{C})$ near $a_0 = \infty$ in the direction of $v_0$ and a path $\gamma$ in $X(\mathbb{C})$ from $v_0$ to $v_0'$. For $1\leq j \leq r$, let $\gamma_j'$ be a path in $X(\mathbb{C})$ from $v_0'$ to $v_j$ and set $\gamma_i := \gamma_i' \cdot \gamma$, where paths are composed from the right. Let $x_j'$ be a small circle around $a_j$ starting from $v_j$ in the opposite clockwise way and set $x_j := \gamma_j^{-1} \cdot x_j' \cdot \gamma_j$. We may assume that paths $x_1' \cdot \gamma_1', \dots , x_{r}' \cdot \gamma_{r}'$ are disjoint each other and that when we make a small circle around $v_0'$ in the opposite clockwise way starting from a point on $\gamma_1'$, we meet successively $\gamma_2', \dots, \gamma_{r}'$.  
 \begin{equation*}
\centering
\includegraphics{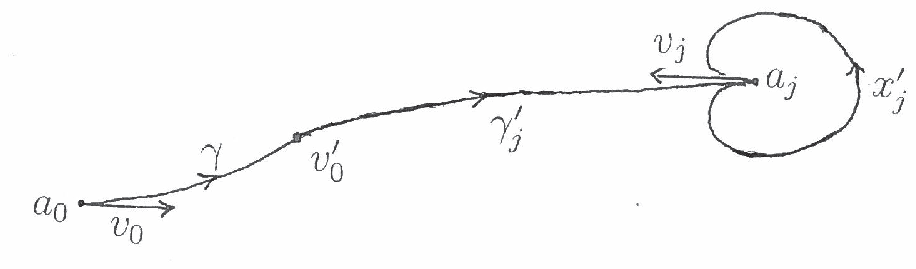}
\end{equation*}
 Then $\pi_1(X(\mathbb{C});v_0)$ is generated by (the homotopy classes of) $x_0, x_1,\dots , x_{r}$ subject to the relation $x_{r} \cdots x_1 x_0= 1$. Hence  $\pi_1(X(\mathbb{C});v_0)$ is identified with the free group $F_r$ generated by $x_1,\dots , x_r$ (the path $x_i$ and the word $x_i$ are identified).  
 For $0\leq  j \leq r$, let $\widehat{\pi_1}^{(l)}(X(\mathbb{C}); v_0, v_j)$ denote the pro-$l$ completion of the set $\pi_1(X(\mathbb{C}); v_0, v_j)$ of homotopy classes of paths in $X(\mathbb{C})$ from $v_0$ to $v_j$. When $v_0 = v_j$, $\widehat{\pi_1}^{(l)}(X(\mathbb{C}); v_0, v_j)$ is denoted by $\widehat{\pi_1}^{(l)}(X(\mathbb{C}); v_0)$ which is  the pro-$l$ completion $\widehat{F_r}^{(l)}$ of $F_r$.  By ([G; XII, Corollaire 5.2]),  $\widehat{\pi_1}^{(l)}(X(\mathbb{C}); v_0)$ is the maximal pro-$l$ quotient of the \'{e}tale fundamental group of $X_{\overline{k}}$ based at $v_0$.
 
Let ${\rm Gal}_k$ denote the absolute Galois group ${\rm Gal}(\overline{k}/k)$ over $k$. The Ihara representation of  ${\rm Gal}_k$ on $\widehat{F_r}^{(l)} = \widehat{\pi_1}^{(l)}(X(\mathbb{C}); v_0)$  is given by the monodromy action as follows. Let $M_A$ be the maximal pro-$l$ extension of $\overline{k}(t)$ unramified outside $A$. Let $t_0 := 1/t$ and $t_j := t - a_j$ for $1\leq j \leq r$. For each $j = 0,\dots , r$, let $\overline{k}\{\{t_j \}\} := \bigcup_{n \geq 1} \overline{k}(( t_j^{1/n}))$ be the field of Puiseux series in $t_j$ with coefficients in $\overline{k}$ and let $\iota_j :  M_A \hookrightarrow \overline{k}\{\{t_j\}\}$ be the natural embedding. Let $M_j$ denote the image of $\iota_j$. For each path $p \in \pi_1(X(\mathbb{C}); v_0, v_j)$ ($0\leq j \leq r$), we have a $\overline{k}(t)$-algebra isomorphism $[p] : M_0 \stackrel{\sim}{\rightarrow} M_j$ by the analytic continuation along $p$. Letting ${\rm Isom}_{\overline{k}(t)}(M_0, M_j)$ denote the set of $\overline{k}(t)$-algebra isomorphisms from $M_0$ to $M_j$, the correspondence $p \mapsto [p]$ induces the bijection
$$ \widehat{\pi_1}^{(l)}(X(\mathbb{C}); v_0, v_j) \stackrel{\sim}{\longrightarrow} {\rm Isom}_{\overline{k}(t)}(M_0, M_j). \leqno{(1.1.1)}$$
 For the particular case that $v_j = v_0$, we have the isomorphism
$$ \widehat{F_r}^{(l)} = \widehat{\pi_1}^{(l)}(X(\mathbb{C}); v_0) \stackrel{\sim}{\longrightarrow} {\rm Gal}(M_0/\overline{k}(t)) \simeq {\rm Gal}(M_A/\overline{k}(t)), \leqno{(1.1.2)}$$
and hence a pro-$l$ word $f \in \widehat{F_r}^{(l)}$ acts on $M_0$ as the monodromy transformations of algebraic functions in $M_A$ along the pro-$l$ path $f$.
The Galois group ${\rm Gal}_k$ acts on $\overline{k}\{\{t_j\}\}$ via the action on Puiseux coefficients. This action stabilizes $M_j$ and so we have a homomorphism
$$ s_j : {\rm Gal}_k \longrightarrow {\rm Gal}(M_j/\overline{k}(t)). $$
Under the identification $(1.1.1)$, we define the action of ${\rm Gal}_k$ on $\widehat{\pi_1}^{(l)}(X(\mathbb{C}); v_0, v_j)$ by
$$ g(p) := s_j(g) \cdot p \cdot s_0(g)^{-1}   $$
for $g \in {\rm Gal}_k$ and $p \in \widehat{\pi_1}^{(l)}(X(\mathbb{C}); v_0, v_j)$. In particular, ${\rm Gal}_k$ acts on $\widehat{F_r }^{(l)} = \widehat{\pi_1}^{(l)}(X(\mathbb{C}); v_0)$ by
$$ g(f) := s_0(g) f  s_0(g)^{-1}  $$
for $g \in {\rm Gal}_k$ and $f \in \widehat{F_r}^{(l)}$ and thus we obtain  the {\em Ihara representation} associated to $A$ and $v_0$
$$ {\rm Ih} = {\rm Ih}_{(A,v_0)} : {\rm Gal}_k \longrightarrow {\rm Aut}(\widehat{F_r}^{(l)}), \leqno{(1.1.3)}$$
where ${\rm Aut}(\widehat{F_r}^{(l)})$ is the group of (topological) automorphisms of $\widehat{F_r}^{(l)}$, which is virtually a pro-$l$ group ([DDMS; Theorem 5.6]). 

For $1\leq j \leq r$, we define the map $f_j : {\rm Gal}_k \rightarrow \widehat{F_r}^{(l)} = \widehat{\pi_1}^{(l)}(X(\mathbb{C}); v_0)$  by
$$ f_j(g) := g(\gamma_j)^{-1} \cdot \gamma_j = s_0(g) \cdot \gamma_j^{-1}\cdot s_j(g)^{-1} \cdot \gamma_j \;\; (g \in {\rm Gal}_k). \leqno{(1.1.4)} $$
It is easy to see that $f_j$ is a $1$-cocycle
$$ f_j(gh) = {\rm Ih}(g)(f_j(h)) f_j(g) \;\; (g, h \in {\rm Gal}_k). \leqno{(1.1.5)}$$
Then the action of ${\rm Gal}_k$ on the generators $x_1, \dots , x_r$ of $\widehat{F_r}^{(l)}$ is given as follows: Let $\chi_l : {\rm Gal}_k \rightarrow \mathbb{Z}_l^{\times}$ be the $l$-cyclotomic character ($\mathbb{Z}_l$ : $l$-adic integers) defined by $g(\zeta_{l^n}) = \zeta_{l^n}^{\chi_l(g)}$ for $g \in {\rm Gal}_k$ and $\zeta_{l^n} := \exp(\frac{2\pi \sqrt{-1}}{l^n})$.\\
\\
{\bf Lemma 1.1.6} ([W1; Proposition 2.2.1]). {\em Notations being as above, we have, for $g \in {\rm Gal}_k$,}
$$ {\rm Ih}(g)(x_0) = x_0^{\chi_l(g)}, \;\; {\rm Ih}(g)(x_j) = f_j(g) x_j^{\chi_l(g)} f_j(g)^{-1} \,(1\leq j \leq r).$$
\\
By Lemma 1.1.6,  the image of the Ihara representation (1.1.3) is in the following pro-$l$ analogue of the pure braid group ([I1])
$$ \begin{array}{l} 
{\cal P}(\widehat{F_r}^{(l)}) \\
:= \left\{ \varphi \in {\rm Aut}(\widehat{F_r}^{(l)}) \, \Big| \begin{array}{l} \varphi(x_j) \sim x_j^{N(\varphi)} \;  (1\leq j \leq r), \varphi(x_1\cdots x_{r}) = (x_1\cdots x_{r})^{N(\varphi)} \\
\;\;  \mbox{for some}\; N(\varphi) \in \mathbb{Z}_l^{\times} \end{array} \right\},
\end{array}$$
where $N \circ {\rm Ih} = \chi_l$.

Let $\Omega_A$ be the subfield of $\overline{k}$ corresponding to the subgroup ${\rm Ker}({\rm Ih})$ of ${\rm Gal}_k$:
$$ \Omega_A := (\overline{k})^{{\rm Ker}({\rm Ih})}, \leqno{(1.1.7)}$$
which we call the {\em Ihara field of definition} for $A$.  It is the smallest field of definition of all finite ramified coverings of $\mathbb{P}^1_{\mathbb{C}}$ unramified outside $A$ whose Galois closures have degree $l$-power (cf. [AI; 3]). Since ${\rm Aut}(\widehat{F_r}^{(l)})$ is virtually a pro-$l$ group, $\Omega_A$ is virtually a pro-$l$ extension of $k$.  Moreover, since ${\rm Ker}({\rm Ih}) \subset {\rm Ker}(\chi_l)$, we have 
$$ k(\zeta_{l^{\infty}}) := \cup_{n\geq 1}k(\zeta_{l^n}) \subset \Omega_A.$$
The ramification in the extension $\Omega_A/k$ was studied by Wojtkowiak ([W4]). We also refer to [AI] for the case that $A$ contains $\{0,1,\infty\}$.  Define the finite set ${\cal S}_A$ of primes of $k$ by
$$ 
{\cal S}_A  :=\left\{ \frak{p} \in {\rm Spm}({\cal O}_k) \; \Big| \begin{array}{l} \; v_{\frak{p}}(l) > 0\; \mbox{or}\;  v_{\frak{p}}(a_i - a_j) > 0 \; \mbox{for some}\; 1\leq i \neq j \leq r \\
\;\; \mbox{or}\; v_{\frak{p}}(a_i) < 0\; \mbox{for some}\; 1\leq i \leq r \end{array} \right\}, \leqno{(1.1.8)} $$
 where $v_{\frak{p}}$ denotes the $\frak{p}$-adic valuation. \\
\\
{\bf Theorem 1.1.9} ([W4; Theorem 7.17]). {\em Notations being as above, the extension $\Omega_A/k$ is unramified outside ${\cal S}_A$.}\\
\\
 {\bf 1.2. The pro-$l$ longitudes of a Galois element.} Let $H$ be the abelianization of $\widehat{F_r}^{(l)}$,  $H :=\widehat{F_r}^{(l)}/[\widehat{F_r}^{(l)}, \widehat{F_r}^{(l)}]$, and let $[f]$ denote  the image of $f \in \widehat{F_r}^{(l)}$ in $H$. We set $X_j := [x_j]$ for $1\leq j \leq r$ so that $H$ is the free $\mathbb{Z}_l$-module with basis $X_1, \dots ,X_r$. For $1\leq j \leq r$, the $j$-th ({\em preferred}) {\em pro-$l$ longitude} of $g \in {\rm Gal}_k$ is defined to be a pro-$l$ word $y_j(g) \in \widehat{F_r}^{(l)}$ which satisfies the following conditions
$$ \left\{ \begin{array}{l}
(1) \; {\rm Ih}(g)(x_j) = y_j(g) x_j^{\chi_l(g)} y_j(g)^{-1}, \\
(2)\;   [y_j(g)] = \prod_{i \neq j} {e_i} X_i \; \mbox{for some} \; e_i \in \mathbb{Z}_l.
\end{array} \right. \leqno{(1.2.1)}
$$
\\
{\bf Lemma 1.2.2} ([KMT; Lemma 3.2.1]). {\em For each $j$ $(1\leq j \leq r)$, the $j$-th pro-$l$ longitude of each Galois element in ${\rm Gal}_k$ exists uniquely. }\\
\\
The following proposition shows that the $j$-th pro-$l$ longitude of $g$ is given by  $f_j(g)$ in (1.1.3). For $z \in k^{\times}$, let $\kappa_z : {\rm Gal}_k \rightarrow \mathbb{Z}_l$ be the Kummer cocycle defined by 
$$g(z^{1/l^n}) = \zeta_{l^n}^{\kappa_z(g)}z^{1/l^n}  \;\; (n \geq 1). \leqno{(1.2.3)}$$ 
We easily see the formula $g^{-1}(z^{1/l^n}) = \zeta_{l^n}^{-\chi_l(g)^{-1}\kappa_z(g)}z^{1/l^n}$. 
\\
\\
{\bf Proposition 1.2.4.} {\em For $1\leq j \leq r$, the pro-$l$ word $f_j(g)$ is the $j$-th pro-$l$ longitude of $g \in {\rm Gal}_k$ and we have }
$$ [f_j(g)] = - \sum_{i \neq j} \kappa_{a_j - a_i}(g) X_i.$$
\\
{\em Proof.} Since the maximal abelian subextension of $M_S$ over $\overline{k}(t)$ is generated by $t_i^{1/l^n}$ for $0\leq i \leq r$ and $n\geq 1$, $[f_j(g)]$  is determined by   its action on $t_i^{1/l^n}$. For $i \neq j$, the monodromy transformation of $t_i^{1/l^n} = (t-a_i)^{1/l^n}$ along $f_j(g) = s_0(g)\cdot \gamma_j^{-1} \cdot s_j(g)^{-1} \cdot \gamma_j$ is given as follows:
$$ \begin{array}{ll}
\iota_0((t-a_i)^{1/l^n}) & \stackrel{\gamma_j}{\longrightarrow} \displaystyle{(a_j-a_i)^{1/l^n} \sum_{m=0}^{\infty} \frac{{1/l^n \choose m} }{(a_j-a_i)^m} (t-a_j)^m} \\
                                    & \stackrel{s_j(g)^{-1}}{\longrightarrow} \displaystyle{ \zeta_{l^n}^{-\chi_l(g)^{-1}\kappa_{a_j-a_i}(g)}\sum_{m=0}^{\infty} \frac{{1/l^n \choose m} }{(a_j-a_i)^m} (t-a_j)^m } \;\; (\mbox{by (1.2.3)})\\
                                     & \stackrel{\gamma_j^{-1}}{\longrightarrow} \zeta_{l^n}^{-\chi_l(g)^{-1}\kappa_{a_j-a_i}(g)} \iota_0((t-a_i)^{1/l^n}) \\
                                     & \stackrel{s_0(g)}{\longrightarrow} \zeta_{l^n}^{-\kappa_{a_j-a_i}(g)} \iota_0((t-a_i)^{1/l^n})\;\;  (\mbox{by}\; \iota_0((t-a_i)^{1/l^n}) \in k\{\{1/t\}\}).\, \\
                                     \end{array}
                                     $$
Similarly, we easily see that $f_j(g)$ acts trivially on $(t-a_j)^{1/l^n}$. Since the monodromy translation of  $(t-a_i)^{1/l^n}$ along $x_j$ is the multiplication by $\zeta_{l^n}$ if $i = j$ and the identity if $i \neq j$, we have 
$$ [f_j(g)] = - \sum_{i \neq j} \kappa_{a_j - a_i}(g) X_i.  \leqno{(1.2.4.1)}$$
By Lemma 1.1.5, (1.2.1), (1.2.4.1) and the uniqueness of the $j$-th pro-$l$ longitude,  $f_j(g)$ is the $j$-th pro-$l$ longitude of $g$ $\;\; \Box$
\\
 
 Let $\Omega_A$ be as in (1.1.7) and let $\sigma \in {\rm Gal}(\Omega_A/k)$. Choosing an extension $g$ of $\sigma$, we set
 $$ f_j(\sigma) := f_j(g) \;\; (1\leq j \leq r). \leqno{(1.2.5)} $$
 \\
 {\bf Proposition 1.2.6.} {\em The definition of $f_j(\sigma)$ in (1.2.5) is independent of the 
 choice of an extension $g \in {\rm Gal}_k$.}\\
 \\
{\em Proof.} Let $g, g' \in  {\rm Gal}_k$ be extensions of $\sigma$. We can write $g' = gh$ for
some $h \in {\rm Gal}(k/\Omega_A)$. By (1.1.5), we have $f_j(g') = {\rm Ih}(g)(f_j(h))f_j(g)$. Since
${\rm Ih}(h) = {\rm id}$, the uniqueness of the pro-$l$ longitude in Lemma 1.2.2 yields
$f_j(h) = 1$ and hence $f_j(g') = f_j(g)$. $\;\; \Box$\\
 \\
 {\bf 1.3. Mod $l$ Milnor invariants of a Galois element.} Let $\widehat{T}$ be the complete tensor algebra of $H$ over $\mathbb{Z}_l$ defined by $\widehat{T} := \prod_{n \geq 0} H^{\otimes n}$, where $H^{\otimes 0} := \mathbb{Z}_l$ and $H^{\otimes n} := H \otimes_{\mathbb{Z}_l} \cdots \otimes_{\mathbb{Z}_l} H$ ($n$ times) for $n \geq 1$. It is nothing but the {\em Magnus algebra}  $\mathbb{Z}_l \langle \langle X_1, \dots , X_r \rangle \rangle$ over $\mathbb{Z}_l$, namely,  the algebra of non-commutative formal power series (called {\em Magnus power series}) over $\mathbb{Z}_l$ with variables $X_1, \dots , X_r$:
$$ \widehat{T} = \prod_{n \geq 0} H^{\otimes n} = \mathbb{Z}_l \langle \langle X_1, \dots , X_r \rangle \rangle.$$
For $n \geq 0$, we set $\widehat{T}(n) := \prod_{m \geq n} H^{\otimes m}$. The {\em degree} of a Magnus power series $\Phi$, denoted by ${\rm deg}(\Phi)$, is defined to be the minimum $n$ such that $\Phi \in \widehat{T}(n)$. We note that  $H^{\otimes n}$ is the free $\mathbb{Z}_l$-module on monomials $X_{i_1}\cdots X_{i_n}$ $(1\leq i_1,\dots , i_n \leq r)$ of degree $n$ and $\widehat{T}(n)$  consists of Magnus power series of degree $\geq n$.  

Let $\mathbb{Z}_l[[\widehat{F_r}^{(l)}]]$ be the complete group algebra of $\widehat{F_r}^{(l)}$ over $\mathbb{Z}_l$ and let $\epsilon_{\mathbb{Z}_l[[\widehat{F_r}^{(l)}]]} : 
\mathbb{Z}_l[[\widehat{F_r}^{(l)}]] \rightarrow \mathbb{Z}_l$ be the augmentation homomorphism with the augmentation ideal $I_{\mathbb{Z}_l[[\widehat{F_r}^{(l)}]]} := {\rm Ker}(\epsilon_{\mathbb{Z}_l[[\widehat{F_r}^{(l)}]]})$.
 The correspondence $x_i \mapsto 1 +X_i$ $(1\leq i \leq r)$ gives rise to the {\em pro-$l$  Magnus isomorphism} of topological $\mathbb{Z}_l$-algebras
$$ \Theta : \mathbb{Z}_l[[\widehat{F_r}^{(l)}]]\; \stackrel{\sim}{\longrightarrow}  \widehat{T} = \mathbb{Z}_l \langle \langle X_1, \dots , X_r \rangle \rangle.  \leqno{(1.3.1)}$$
Here $(I_{\mathbb{Z}_l[[\widehat{F_r}^{(l)}]]})^n$ corresponds, under $\Theta$,  to $\widehat{T}(n)$ for $n \geq 0$. For $\alpha \in \mathbb{Z}_l[[\widehat{F_r}^{(l)}]]$, $\Theta(\alpha)$ is called the {\it pro-$l$ Magnus expansion} of $\alpha$. In the following, for a multi-index $I = (i_1 \cdots i_n)$,  $1\leq i_1,\dots , i_n \leq r$, we set
$$ |I| := n \; \mbox{and} \; X_I  := X_{i_1} \cdots X_{i_n}.$$
We call the coefficient of $X_I$ in $\Theta(\alpha)$ the {\em $l$-adic Magnus coefficient} of $\alpha$ for $I$  and denote it by $\mu^{(l)}(I;\alpha)$. So we have
$$ \Theta(\alpha) = \epsilon_{\mathbb{Z}_l[[\widehat{F_r}^{(l)}]]}(\alpha) + \sum_{|I|\geq 1} \mu^{(l)}(I;\alpha) X_I.  $$
Taking mod $l$ in (1.3.1), we have the {\em mod $l$ Magnus isomorphism}
$$ \Theta_l : \mathbb{F}_l[[\widehat{F_r}^{(l)}]]\; \stackrel{\sim}{\longrightarrow}  \widehat{T}\otimes_{\mathbb{Z}_l}\mathbb{F}_l  = \mathbb{F}_l \langle \langle X_1, \dots , X_r \rangle \rangle \leqno{(1.3.2)} $$
so that for $\alpha \in \mathbb{F}_l[[\widehat{F_r}^{(l)}]]$, we have 
$$ \Theta_l(\alpha) = \epsilon_{\mathbb{F}_l[[\widehat{F_r}^{(l)}]]}(\alpha) + \sum_{|I|\geq 1} \mu_l(I;\alpha) X_I, $$
where $\epsilon_{\mathbb{F}_l[[\widehat{F_r}^{(l)}]} : \mathbb{F}_l[[\widehat{F_r}^{(l)}]] \rightarrow \mathbb{F}_l$ is the augmentation homomorphism and $\mu_l(I;\alpha) := \mu
(I;\alpha)$ mod $l$. Let $\{ \widehat{F_r}^{(l)}(d) \}_{d\geq 1}$ be the Zassenhaus filtration of $ \widehat{F_r}^{(l)}$ defined by $\widehat{F_r}^{(l)}(d) := \widehat{F_r}^{(l)} \cap 1 + (I_{\mathbb{F}_l[[\widehat{F_r}^{(l)}]})^d$, where $I_{\mathbb{F}_l[[\widehat{F_r}^{(l)}]]} := {\rm Ker}(\epsilon_{\mathbb{F}_l[[\frak{F}_r]]})$. For $f \in  \widehat{F_r}^{(l)}$, we have
$$ f \in \widehat{F_r}^{(l)}(d)  \Longleftrightarrow \mu_l(I;f) = 0 \; \mbox{for} \; |I| < d \;\; \mbox{i.e.,}\, {\rm deg}(\Theta_l(f-1)) \geq d. $$
The following inductive formula for $\widehat{F_r}^{(l)}(d)$ is known ([DDMS, 12.9]):
$$ \widehat{F_r}^{(l)}(d) = (\widehat{F_r}^{(l)}([d/l]))^l \prod_{i+j=d} [\widehat{F_r}^{(l)}(i),\widehat{F_r}^{(l)}(j)], \leqno{(1.3.3)}
$$
where $[d/l]$ stands for the least integer $m$ such that $ml \geq d$. \\

Now, following the case for pure braids ([MK; Chapter 6, 4], [Kd; Chapter 1]), we will define the $l$-adic Milnor numbers of $g \in {\rm Gal}_{k}$ by the $l$-adic Magnus coefficients of the $i$-th longitude $y_i(g)$: Let $I = (i_1\cdots i_n)$ be a multi-index, where  $1\leq i_1, \dots , i_n \leq r$ and $|I| = n \geq 1$. The {\it $l$-adic Milnor number} of $g \in {\rm Gal}_{k}$ for $I$, denoted by $\mu^{(l)}(g; I) = \mu^{(l)}(g; i_1\cdots i_n)$, is defined by the $l$-adic Magnus coefficient of the pro-$l$ longitude $f_{i_n}(g)$ for $I' := (i_1 \cdots i_{n-1})$:
$$  \mu^{(l)}(g; I) := \mu^{(l)}(I'; f_{i_n}(g)).  $$
Here we set $\mu^{(l)}(g; I) := 0$ if $|I| = 1$.  In this paper, we shall use {\em mod $l$ Milnor number} $\mu_l(g;I)$ of $g \in {\rm Gal}_{k}$ for $I$, which is  defined by
$$ \mu_l(g;I) := \mu_l(I';f_{i_n}(g)) := \mu^{(l)}(g;I) \; \mbox{mod}\; l. \leqno{(1.3.4)}$$
By the proof of [KMT; Theorem 3.2.8], we have the following\\
\\
{\bf Theorem 1.3.5.} {\em  Let $g, h \in {\rm Gal}_k$ satisfying $\chi_l(g) \equiv \chi_l(h) \equiv 1$ mod $l$. Let $I = (i_1\cdots i_n)$ be a multi-index. We assume that $\mu_l(g;J) = 0$ for any $J = (j_1\cdots j_m)$ with $\{j_1,\dots, j_m \} \subsetneqq \{ i_1,\dots, i_n \}$. Then we have}
$$\mu_l(hgh^{-1};I) = \mu_l(g;I). $$
When the conditions in Theorem 1.3.5 are satisfied, we call $\mu_l(g;I)$ the {\em mod $l$ Milnor invariant} of $g$ for $I$.

Let $\Omega_A$  be the Ihara field of definition for $A$ in (1.1.7). By Proposition 1.2.6, mod $l$ Milnor number $\mu _l(\sigma; I)$ of $\sigma \in {\rm Gal}(\Omega_A/k)$ for a multi-index $I$ is
well defined by $\mu_l(g; I)$ for an extension $g \in {\rm Gal}_k$ of $\sigma$. Let ${\cal S}_A$ be as in (1.1.8). Let $\frak{p} \in {\rm Spm}({\cal O}_k) \setminus {\cal S}_A$ and let $\frak{P}$ be an extension of $\frak{p}$ to $\Omega_A$. Since $\frak{P}$ is unramified in $\Omega_A/k$ by Theorem 1.1.9, we have  the Frobenius automorphism $\sigma_{\frak{P}} \in {\rm Gal}(\Omega_A/k)$  of $\frak{P}$ over $k$. We then have mod $l$ Milnor number $\mu_l(\sigma_{\frak{P}};I)$ for a multi-index $I$.\\
\\
{\bf Corollary 1.3.6.} {\em Notations being as above, suppose ${\rm N}\frak{p} \equiv 1 \; \mbox{mod}\; l$. Let $I = (i_1\cdots i_n)$ be a multi-index. We assume that $\mu_l(\sigma_{\frak{P}};J) = 0$ for any $J = (j_1\cdots j_m)$ with $\{j_1,\dots, j_m \} \subsetneqq \{ i_1,\dots, i_n \}$. Then $\mu_l(\sigma_{\frak{P}};I)$ is independent of the choice of an extension $\frak{P}$ and hence it is denoted by  $\mu_l(\sigma_{\frak{p}};I)$}\\
\\
{\em Proof.} This follows from Theorem 1.3.5 and $\chi_l(\sigma_{\frak{P}}) = {\rm N}\frak{p} \equiv 1 \; \mbox{mod}\; l.$ $\;\; \Box$\\
\\
{\bf Theorem 1.3.7.} {\em Notations being as above, for $1\leq i \leq r$, we have}
$$ \mu_l(\sigma_{\frak{p}}; ii) = 0.$$
{\em For $1\leq i \neq j \leq r$, we have }
$$ k((a_j-a_i)^{1/l^n}) \subset \Omega_A\;\; (n \geq 1) $$
and 
$$ \displaystyle{ \zeta_l^{\mu_l(\sigma_{\frak{p}}; ij)} = \left( \frac{a_j- a_i}{\frak{p}} \right)_l^{-1}} $$
{\em for $\frak{p} \notin {\cal S}_A$ with ${\rm N}\frak{p} \equiv 1 \; \mbox{mod} \; l$. Here $\left( \frac{\cdot}{\frak{p}} \right)_l$ denotes the $l$-th power residue symbol in $k_{\frak{p}}$.}\\
\\
{\em Proof.} The first assertion follows from Proposition 1.2.4. For the second assertion, it suffices to show that $g((a_j-a_i)^{1/l^n}) = (a_j-a_i)^{1/l^n}$ for any $g \in {\rm Gal}(\overline{k}/\Omega_A)$. Since ${\rm Gal}(\overline{k}/\Omega_A) = {\rm Ker}({\rm Ih})$ by (1.1.7), we have $f_j(g) = 1$ for $1\leq j \leq r$ by (1.2.1), Lemma 1.2.2 and Proposition 1.2.4. By Proposition 1.2.4 again,
we have $\kappa_{a_j-a_i}(g) = 0$ for $i \neq j$. By (1.2.3), $g((a_j-a_i)^{1/l^n}) = (a_j-a_i)^{1/l^n}$. 

We note by the second assertion and Theorem 1.1.9 that $k_{\frak{p}}((a_i-a_j)^{1/l})$ is an unramified extension of $k_{\frak{p}}$ for $\frak{p} \notin {\cal S}_A$. By (1.2.3), Proposition 1.2.4 and (1.3.4), the third assertion is obtained as follows:
$$ \begin{array}{ll}
\zeta_l^{\mu_l(\sigma_{\frak{p}}; ij)} & = \zeta_l^{\mu_l(i;f_j(\sigma_{\frak{P}}))}\\
                                                             & = \zeta_l^{- \kappa_{a_j-a_i}(\sigma_{\frak{P}})}\\
                                                            & = \displaystyle{ \left( \frac{\sigma_{\frak{P}}((a_j-a_i)^{1/l})}{(a_j-a_i)^{1/l}} \right)^{-1}  }\\
                                                            & = \displaystyle{ \left( \frac{a_j- a_i}{\frak{p}} \right)_l^{-1} } \;\; \Box
                                                            \end{array}
$$
\vspace{0.5cm}\\
{\bf Remark 1.3.8.} (1) By the relation between Magnus coefficients and Massey products ([Dw],[St]), it was shown in [KMT; $\S 3.3$] that the mod $l$ Milnor invariants of a Galois element $g$ are expressed by Massey products in the mod $l$ cohomology of the pro-$l$ link group of $g$ defined by
$$ \Pi_A(g)  := \langle x_1, \dots , x_r \, |  \,  x_1^{1-\chi_l(g)}[x_1^{-1},f_1(g)^{-1}] = \cdots = x_r^{1-\chi_l(g)}[x_r^{-1},f_r(g)^{-1}] = 1 \rangle. $$
(2) Let $E : \widehat{F_r}^{(l)} \rightarrow \mathbb{Q}_l \langle \langle X_1,\dots , X_r \rangle \rangle$ be the embedding defined by $E(x_i) := \exp(X_i) = 1 + X_i + \frac{1}{2}X_i^2 + \cdots + \frac{1}{n !} X_i^n + \cdots$.  In a series of papers [NW], [W1] $\sim$ [W4], Wojtkowiak has studied the coefficients of Lie elements in the series  $\log E(f_j(g)^{-1})$, called the {\em $l$-adic iterated integrals}. Our $l$-adic Milnor numbers are expressed by $l$-adic iterated integrals, and $l$-adic iterated integrals, vice versa. \\
 
 \begin{center}
 {\bf 2. Dilogarithmic mod $l$ Heisenberg ramified coverings of $\mathbb{P}^1$}
 \end{center}
 
 In this section, we introduce certain mod $l$ Heisenberg extensions of $k(t)$, called the dilogarithmid mod $l$ Heisenberg extensions, which will be used later in the section 4. We explain the analogies between our mod $l$ Heisenberg coverings and the dilogarithmic function from cohomological viewpoint. We assume that the number field $k$ contains $\zeta_l = \exp(\frac{2\pi \sqrt{-1}}{l})$. \\
\\
{\bf 2.1. Mod $l$ Heisenberg branched coverings of $\mathbb{P}^1$.} Let $k(t)$ be the function field of the projective $t$-line $\mathbb{P}^1$ over $k$. For $c \in k^{\times}$, let ${\cal K}^{(l)} := {\cal K}_{\{c,t\}}^{(l)}$ be the extension of $k(t)$ defined by
 $$ {\cal K}^{(l)} := {\cal K}_{\{c,t\}}^{(l)} := k(t)(t^{1/l}, (c^l-t)^{1/l}). \leqno{(2.1.1)}$$
It is a Kummer extension of $k(t)$ such that the Galois group ${\rm Gal}({\cal K}^{(l)}/k(t))$ is isomorphic to $\mathbb{Z}/l\mathbb{Z} \times \mathbb{Z}/l\mathbb{Z}$ generated by $\alpha, \beta$ defined by
$$ \begin{array}{ll}
\alpha(t^{1/l}) := \zeta_l \,t^{1/l}, & \alpha((c^l-t)^{1/l}) := (c^l-t)^{1/l},\\
\beta(t^{1/l}) :=t^{1/l}, & \beta((c^l-t)^{1/l}) :=  \zeta_l \,(c^l-t)^{1/l},\\
\end{array}$$
and it is unramified outside $ \infty, 0$ and $c^l$. The ramification index of these points are $l$. A non-singular projective curve ${\cal C}^{(l)}$ over $k$ with function field ${\cal K}^{(l)}$ is given by the Fermat plane curve
$$ \; X^l + Y^l = (cZ)^l$$
in $\mathbb{P}^2$ and the covering map ${\cal C}^{(l)} \rightarrow \mathbb{P}^1$ is given by
$$ (X:Y:Z) \, \mapsto \, (X^l:Z^l).$$

We set 
$$ \varepsilon_l(t) := \prod_{i=1}^{l-1} (c -\zeta_l^i t^{1/l})^i $$
and define the extension ${\cal R}^{(l)} := {\cal R}_{\{c,t\}}^{(l)}$ of ${\cal K}^{(l)}$  by
$$ {\cal R}^{(l)} := {\cal R}_{\{c,t\}}^{(l)} := {\cal K}^{(l)}(\varepsilon_l(t)^{1/l}) = k(t)(t^{1/l}, (c^l-t)^{1/l}, \varepsilon_l(t)^{1/l}). \leqno{(2.1.2)}$$
It is a cyclic Kummer extension of ${\cal K}^{(l)}$ of degree $l$ whose Galois group  ${\rm Gal}({\cal R}^{(l)}/{\cal K}^{(l)})$ is generated by $\delta$ defined by
$$ \delta(\varepsilon_l(t)^{1/l}) := \zeta_l\, \varepsilon_l(t)^{1/l},$$
and in which  only primes $ (c -\zeta_l^i t^{1/l})$ $(0\leq i \leq l-1)$ of ${\cal K}^{(l)}$, which are all lying over $t=c^l$, can be ramified in ${\cal R}^{(l)}$. Let ${\cal D}^{(l)}$ be a non-singular projective curve whose function field is ${\cal R}^{(l)}$. For $l = 2$ and $3$, concrete defining equations for ${\cal D}^{(l)}$ are given as follows.\\
\\
{\bf Example 2.1.3.} Let $l = 2$. By setting $U/W = \sqrt{c - \sqrt{t}}$ and $V/W = \sqrt{c + \sqrt{t}}$, we can take a non-singular projective model ${\cal D}^{(2)}$ of ${\cal R}^{(2)}$ by the plane curve
$$ U^2 + V^2 = 2c W^2$$
and hence the genus of ${\cal D}^{(2)}$ is $0$. The covering map ${\cal D}^{(2)} \rightarrow {\cal C}^{(2)}$ is given by
$$ (U:V:W) \, \mapsto \, (cW^2-U^2: UV : W^2),$$
which is ramified at $(0:\pm \sqrt{2c} : 1) \in {\cal C}^{(2)}$. 

Let $l = 3$. By setting $U = \sqrt[3]{c - \zeta_3 \sqrt[3]{t}}, V = \sqrt[3]{c - \zeta_3^2 \sqrt[3]{t}}$ and $W = \sqrt[3]{c - \sqrt[3]{t}}$, we can take a non-singular projective model ${\cal D}^{(3)}$ of ${\cal R}^{(3)}$ by the plane curve
$$ \zeta_3^2 U^3 + V^3 = - \zeta_3W^3$$
and hence the genus of ${\cal D}^{(3)}$ is $1$. The covering map ${\cal D}^{(3)} \rightarrow {\cal C}^{(3)}$ is given by
$$ (U:V:W) \, \mapsto \, (c(U^3-W^3): c(1-\zeta_3)UVW : U^3- \zeta_3 W^3),$$
which is unramified. 
\\
\\
{\bf Theorem 2.1.4.} {\em Notations being as above, ${\cal R}^{(l)}$ is a Galois extension of $k(t)$ such that Galois group ${\rm Gal}({\cal R}^{(l)}/k(t))$ is isomorphic to the mod $l$ Heisenberg group}
$$ H(\mathbb{F}_l) := \left\{ \left( \begin{array}{ccc} 1 & * & * \\ 0 & 1 & * \\ 0 & 0 & 1 \end{array} \right) \Big| \; * \in \mathbb{F}_l \; \right\},$$
 {\em and it is unramified outside $\infty, 0$ and $c^l$.}\\
 \\
 {\em Proof.} The assertion about the ramification in ${\cal R}^{(l)}/k(t)$ follows immediately from those in ${\cal K}^{(l)}/k(t)$ and ${\cal R}^{(l)}/{\cal K}^{(l)}$. For $0\leq j < l$, we see that
 $$ \displaystyle{ \alpha^j(\varepsilon_l(t)) = \frac{\{ (c-t^{1/l})\cdots (c-\zeta^{j-1}t^{1/l})\}^l}{(c^l -t)^j}\, \varepsilon_l(t), \;\; \beta^j(\varepsilon_l(t)) = \varepsilon_l(t),}$$
 from which any congugate of $\varepsilon(t)^{1/l}$ over $k(t)$ lies in ${\cal R}^{(l)}$ and so ${\cal R}^{(l)}$ is a Galois extension of $k(t)$. 
 We define the extensions $\tilde{\alpha}, \tilde{\beta} \in {\rm Gal}({\cal R}^{(l)}/k(t))$ of $\alpha, \beta \in {\rm Gal}({\cal K}^{(l)}/k(t))$, respectively, by
 $$ \begin{array}{lll}
   \tilde{\alpha}(t^{1/l}) := \zeta_l \, t^{1/l}, & \tilde{\alpha}((c^l -t)^{1/l}) := (c^l -t)^{1/l}, & \displaystyle{\tilde{\alpha}(\varepsilon_l(t)^{1/l}) := \frac{c-t^{1/l}}{(c^l -t)^{1/l}} \varepsilon_l(t)^{1/l}, } \\
  \tilde{\beta}(t^{1/l}) := t^{1/l}, & \tilde{\beta}((c^l -t)^{1/l}) :=  \zeta_l \, (c^l -t)^{1/l}, & \tilde{\beta}(\varepsilon_l(t)^{1/l}) := \varepsilon_l(t)^{1/l},
 \end{array}
 \leqno{(2.1.4.1)}
 $$
where we easily verify that $\tilde{\alpha}^l = \tilde{\beta}^l = 1$. By the straightforward computation, we have
 $$ [\tilde{\alpha}, \tilde{\beta}](t^{1/l}) = t^{1/l},  \, [\tilde{\alpha},\tilde{\beta}]((c^l-t)^{1/l}) = (c^l-t)^{1/l},  [\tilde{\alpha},\tilde{\beta}](\varepsilon_l(t)^{1/l}) = \zeta_l \varepsilon_l(t)^{1/l} \leqno{(2.1.4.2)}$$
and so $[\tilde{\alpha}, \tilde{\beta}] = \delta$. Therefore the correspondence
$$ \tilde{\alpha} \mapsto \left( \begin{array}{ccc} 1 & 1 & 0 \\ 0 & 1 & 0 \\ 0 & 0 & 1 \end{array}\right), \; \tilde{\beta} \mapsto  \left( \begin{array}{ccc} 1 & 0 & 0 \\ 0 & 1 & 1\\ 0 & 0 & 1 \end{array}\right), $$
induces the isomorphism
$$ {\rm Gal}({\cal R}^{(l)}/k(t)) \stackrel{\sim}{\longrightarrow} H(\mathbb{F}_l).  \;\; \Box $$
 \\
 We call the extension ${\cal R}^{(l)}/k(t)$ the {\em dilogarithmic mod $l$ Heisenberg extension}, and call the ramified covering ${\cal D}^{(l)} \rightarrow \mathbb{P}^1_k$ (resp. the (unramified) covering ${\cal D}^{(l)}|_X \rightarrow X$ for $X := \mathbb{P}^1 \setminus \{ \infty, 0, c^l \}$) the {\em dilogarithmic mod $l$ Heisenberg ramified covering} (resp. the {\em dilogarithmic mod $l$ Heisenberg covering}). A mod $l$ Heisenberg extension (resp. (ramified) covering) will also be called simply an $H(\mathbb{F}_l)$-extension (resp. $H(\mathbb{F}_l)$-(ramified) covering). The reason why we call ``{\em dilogarithmic}" will be explained in the next subsection 2.2.
 \\

By Theorem 2.1.4, we have the surjective homomorphism 
$${\rm Gal}(M_A/k(t)) \rightarrow {\rm Gal}({\cal R}^{(l)}/k(t)),$$
when $A$ contains $\infty, 0$ and $c^l$. Composing with it the natural homomorphism $\widehat{\pi_1}^{(l)}(\mathbb{P}^1(\mathbb{C}) \setminus A) \rightarrow  {\rm Gal}(M_A/k(t))$ obtained by the isomorphism (1.1.2) and the inclusion ${\rm Gal}(M_A/\overline{k}(t)) \subset {\rm Gal}(M_A/k(t))$, we have the homomorphism
$$\rho : \widehat{\pi_1}^{(l)}(\mathbb{P}^1(\mathbb{C}) \setminus A; v_0) \longrightarrow {\rm Gal}({\cal R}^{(l)}/k(t)).$$
Let $a_0 := \infty, a_1 = 0$ and $a_2 = c^l$ and let $x_1$ and $x_2$ be the loops around $0$ and $c^l$, respectively, as in Subsection 1.1. Then we have
 $$ \rho(x_1) = \tilde{\alpha}, \; \rho(x_2) = \tilde{\beta}. \leqno{(2.1.5)}$$
 \\
 {\bf Corollary 2.1.6.} {\em The monodromy transformation of $\varepsilon_l(t)^{1/l}$ along the pro-$l$ path $[x_i,x_j]$ $(resp.\,  x_i^l)$ $(1\leq i < j \leq r)$ is given by}
 $$ \varepsilon_l(t)^{1/l} \mapsto \left\{ \begin{array}{ll} \zeta_l  \varepsilon_l(t)^{1/l} & \; (i,j) = (1,2),  \\ \varepsilon_l(t)^{1/l} & \; otherwise, \end{array}   \right. \;\; (resp. \; \varepsilon_l(t)^{1/l} \mapsto \varepsilon_l(t)^{1/l}). $$
 \\
 {\em Proof.} The assertion for  the monodromy along  $[x_i, x_j]$ follows from (2.1.4.2) and  (2.1.5) when $(i,j) = (1,2)$, and from $\rho([x_i,x_j]) = {\rm id}$ when $(i,j) \neq (1,2)$. The assertion for the monodromy along  $x_i^l$ follows from (2.1.4.1) and $\rho(x_i) = {\rm id}$ when $i \neq 1,2$.  $\;\; \Box$\\
 \\
 {\bf Theorem 2.1.7.} {\em Let $a \in k \setminus \{ 0, c^l \}$ and let  $A := \{ a_0 = \infty, a_1 := 0, a_2 := c^l, a_3 := a\}$. Let  $\Omega_A$ be the Ihara field of definition for $A$ in (1.1.7).  Then we have}
 $$ {\cal R}^{(l)}_{\{c,a\}} \subset \Omega_A.$$
 \\
 {\em Proof.}  It suffices to show that
$$ g(a^{1/l}) = a^l, \; g((c^l - a)^{1/l}) = (c^l - a)^{1/l} \;\; \mbox{and} \;\; g(\varepsilon_l(a)^{1/l}) = \varepsilon_l(a)^{1/l} \leqno{(2.1.7.1)}$$
for any $g \in {\rm Gal}(\overline{k}/\Omega_A)$. Since ${\rm Gal}(\overline{k}/\Omega_A) = {\rm Ker}({\rm Ih})$, we have $f_j(g) = 1$ for $1\leq j \leq 3$ by (1.2.1), Lemma 1.2.2 and Proposition 1.2.4. By Proposition 1.2.4, noting $a_3 - a_1 = a$ and $a_2 - a_3 = c^l - a$, we have $\kappa_a(g) = \kappa_{c^l -a}(g) = 1$, which yields the first 2 equalities in (2.1.7.1). To prove the 3rd equality in (2.1.7.1), we first note that for each $n \geq 0$, there is $\Phi_n(X_1,X_2) \in k(X_1,X_2)$ such that 
$$ \frac{d^n}{dt^n} \varepsilon_l(t)^{1/l}\Big|_{t=a} = \varepsilon_l(a)^{1/l} \Phi_n(\varepsilon_l(a), a^{1/l}).$$
Using this and the 1st equality in (2.1.7.1), the monodromy transformation of $\varepsilon_l(t)^{1/l}$ along  $f_3(g) =  s_0(g) \cdot \gamma_3^{-1}\cdot s_3(g)^{-1} \cdot \gamma_3$ 
is computed as follows:
$$ \begin{array}{ll}
\iota_0(\varepsilon_l(t)^{1/l}) & \stackrel{\gamma_3}{\longrightarrow} \displaystyle{ \varepsilon_l(a)^{1/l} \sum_{n=0}^{\infty} \Phi_n(\varepsilon_l(a), a^{1/l}) (t-a)^n} \\
                                    & \stackrel{g^{-1}}{\longrightarrow} \displaystyle{ \frac{g^{-1}(\varepsilon_l(a)^{1/l})}{ \varepsilon_l(a)^{1/l}}  \varepsilon_l(a)^{1/l} \sum_{n=0}^{\infty}  \Phi_n(\varepsilon_l(a), a^{1/l}) (t-a)^n}  \\
                                     & \stackrel{\gamma_3^{-1}}{\longrightarrow} \displaystyle{ \frac{g^{-1}(\varepsilon_l(a)^{1/l})}{ \varepsilon_l(a)^{1/l}}   \iota_0(\varepsilon_l(t)^{1/l}) }\\
                                     & \stackrel{s_0(g)}{\longrightarrow}  \displaystyle{  \frac{\varepsilon_l(a)^{1/l}}{ g(\varepsilon_l(a)^{1/l})}  \iota_0(\varepsilon_l(t)^{1/l}).} \;\; \Box\\
                                     \end{array}
                                     $$ 
Since $f_3(g) =1$, we have $g(\varepsilon_l(a)^{1/l}) = \varepsilon_l(a)^{1/l}$ for any $g \in {\rm Gal}(\overline{k}/\Omega_A)$. $\;\; \Box$
\\
 \\
 {\bf Remark 2.1.8.} (1) The dilogarithmic $H(\mathbb{F}_l)$-extension ${\cal R}^{(l)}$ of $k(t)$ is a special case of Anderson-Ihara's elementary extensions ([AI]) and Wojtkowiak's polylogarithmic extensions ([W5; 3]).\\
(2)For the case that $A$ contains $\infty, 0$ and $1$, it was shown in [AI] that $\Omega_A$ is generated over $k$ by algebraic numbers generalizing  higher circular $l$-units.\\
 \\
 {\bf 2.2. Gerbes and analogies with the dilogarithmic function.} In this subsection, we explain the reason why we call ${\cal D}^{(l)}|_X \rightarrow X$ the {\em dilogarithmic} $H(\mathbb{F}_l)$-covering. It comes from some analogies with the dilogarithmic function, which also explain a geometric meaning of our $H(\mathbb{F}_l)$-coverings. The analogies we discuss in this subsection were suggested by Brylinski's work ([Br1], [Br2]), and we refer to [Br1] for materials on Deligne cohomology and gerbes. 
 
  First, let us recall the dilogarithmic function side. Let $f_1$ and $f_2$ be invertible holomorphic functions on $X := \mathbb{P}^1(\mathbb{C}) \setminus A$. Let $H^n_D(X, \mathbb{Z}(n))$ ($n\geq 1$) denote the holomorphic Deligne cohomology, the $n$-th hypercohomology of the Deligne complex $(2\pi \sqrt{-1})^n\mathbb{Z} \rightarrow {\cal O}_X \stackrel{d}{\rightarrow} \cdots \stackrel{d}{\rightarrow} \Omega^{n-1}_X$ ([Br1; Definition 1.5.9]). Since $H^0(X, {\cal O}_X^{\times}) \simeq H^1_D(X,\mathbb{Z}(1))$ ([ibid; Proposition 1.5.10]), 
 each $f_i$ defines a class $c(f_i) \in H^1_D(X,\mathbb{Z}(1))$. Recall that $H^2_D(X,\mathbb{Z}(2))$ classifies isomorphism classes of holomorphic line bundles over $X$ with holomorphic connection ([ibid; Theorem 2.2.20])
$$ H^2_D(X,\mathbb{Z}(2)) \; \simeq \; \left \{\begin{array}{l} \mbox{isom. classes of holomorphic line bundles }\\
\mbox{over $X$ with holomorphic connection} \end{array}\right \}.$$
Hence the cup product $c(f_1) \cup c(f_2)$ defines an isomorphism class of holomorphic line bundle with holomorphic connection, denoted by $(f_1, f_2)$,  which we call the Deligne line bundle. In more concrete terms, the transition function of $(f_1,f_2)$ is given by $f_2^{\log_i f_1 - \log_j f_1}$ on $U_i \cap U_j$ and the connection $1$-form is given by $\log_i f d\log f_2$ on $U_i$, where $X = \bigsqcup_i U_i$ is an open cover and $\log_i f$ is a chosen branch of $\log f$ on $U_i$. The map $\{ f_1, f_2 \} \rightarrow (f_1,f_2)$ is known to be  the Bloch-Beilinson regulator ([Be], [Bl; $\S 1$])
$$K_2(X) \longrightarrow H^2_D(X,\mathbb{Z}(2)).$$ 
 We note that $(f_1, f_2) = 1$ if and only if there is a trivialization of $(f_1, f_2)$, namely, a horizontal section.
In particular, 
let $A = \{ \infty, 0, 1\}$ and $f_1 = 1-t, f_2 = t$. Then the dilogarithmic function
$$ Li_2(t) = -\int_0^t \log(1-t) d\log t = \sum_{n=1}^{\infty} \frac{t^n}{n^2}$$
gives a horizontal section of $(1-t,t)$ ([Bl; $\S 1$], [De; Example 3.5]). The triviality $(1-t,t)=1$ reflects the Steinberg relation in $K_2(X)$.

Next, let us see the Heisenberg covering side.  Let $f_1$ and $f_2$ be invertible regular functions on $X = \mathbb{P}_k^1 \setminus A$. Let $H^n_{\scriptsize \mbox{\'{e}t}}(X,\mathbb{F}_l)$ denote the $n$-th \'{e}tale cohomology group. Since $k$ contains $\zeta_l$, we note $H^n_{\scriptsize \mbox{\'{e}t}}(X, \mu_l^{\otimes i}) = H^n_{\scriptsize \mbox{\'{e}t}}(X,\mathbb{F}_l)$, where $\mu_l$ is the \'{e}tale sheaf of $l$-th roots of unity on $X$. By Kummer class map $H^0_{\scriptsize \mbox{\'{e}t}}(X, \mathbb{G}_m) \rightarrow H^1_{\scriptsize \mbox{\'{e}t}}(X,\mathbb{F}_l)$, each $f_i$ defines a class $c(f_i) \in H^1_{\scriptsize \mbox{\'{e}t}}(X,\mathbb{F}_l)$. Recall that $H^2_{\scriptsize \mbox{\'{e}t}}(X,\mathbb{F}_l)$ classifies equivalence classes of gerbes over $X$ with band $\mathbb{F}_l$ ([Br1; Theorem 5.2.8], [Gi])
$$ H^2_{\scriptsize \mbox{\'{e}t}}(X,\mathbb{F}_l) \; \simeq \; \left\{\mbox{equiv. classes of gerbes over $X$ with band $\mathbb{F}_l$} \right\}.$$
Hence the cup product $c(f_1) \cup c(f_2)$ defines an isomorphism class  of gerbes with band $\mathbb{F}_l$, denoted by $(f_1, f_2)_l$. 
In more concrete terms, $(f_1,f_2)_l$ is the gerbe associated to the central extension of group schemes over $X$
$$ 1 \longrightarrow \mathbb{F}_l \longrightarrow H(\mathbb{F}_l) \longrightarrow \mathbb{F}_l^{\oplus 2}\longrightarrow 1$$
and the $\mathbb{F}_l^{\oplus 2}$-covering ${\cal C}^{(l)}(f_1,f_2) := {\rm Spec} (k[t, (t-a_i)^{-1} (1\leq i \leq r), f_1^{1/l}, f_2^{1/l}]) \rightarrow X.$ 
So the gerbe $(f_l, f_2)_l$ is the obstruction to lifting of  the $\mathbb{F}_l^{\oplus 2}$-covering ${\cal C}^{(l)}(f_1,f_2) \rightarrow X$
to an $H(\mathbb{F}_l)$-covering ([Br1; 5.2], [Br2; 5]). The map $\{ f_1, f_2 \} \rightarrow (f_1,f_2)_l$ is known to be the Soul\'{e} regulator ([So])
$$K_2(X)\otimes \mathbb{F}_l \longrightarrow H^2_{\scriptsize \mbox{\'{e}t}}(X,\mathbb{F}_l).$$ 
We note that $(f_1, f_2) = 1$ if and only if there is a trivialization of $(f_1, f_2)_l$, namely, an $H(\mathbb{F}_l)$-covering over $X$, which lifts ${\cal C}^{(l)}(f_1,f_2) \rightarrow X$.
Without loss of generality, we may assume $c = 1$ for ${\cal C}^{(l)}$ and ${\cal D}^{(l)}$, and let $A = \{ \infty, 0, 1\}$ and $f_1 = 1 -t, f_2 = t$. Then  ${\cal C}^{(l)}(1 -t,t) = {\cal C}^{(l)}$ and so the $H(\mathbb{F}_l)$-covering ${\cal D}^{(l)}|_X \rightarrow X$ gives a trivialization of $(1-t,t)_l$. 
 
Summing up, we have the following comparison. So our $H(\mathbb{F}_l)$-covering ${\cal D}^{(l)}$ over $X$  may be regarded as a categorical higher analog of the dilogarithmic function.
\\
 \begin{center}
 \begin{tabular}{ | c | c |} 
\hline
Deligne line bundle &  Gerbe associated to $H_3(\mathbb{F}_l)$ \\
$(f_1, f_2) \in H^2_D(X, \mathbb{Z}(2))$ &  $(f_1,f_2)_l \in H^2_{\scriptsize \mbox{\'{e}t}}(X, \mathbb{F}_l)$\\   
\hline
Bloch-Beilinson regulator &  Soul\'{e} regulator \\
$K_2(X) \rightarrow  H^2_D(X, \mathbb{Z}(2))$ & $K_2(X) \rightarrow H^2_{\scriptsize \mbox{\'{e}t}}(X, \mathbb{F}_l)$ \\
\hline
Trivialization of $(1-t, t)$: &   Trivialization of $(1-t, t)_l$: \\ 
  Dilogarithmic function & Mod $l$ Heisenberg covering\\
  $Li_2(t) = -\int_0^t \log(1-t) \frac{dt}{t}$ & ${\cal D}^{(l)}|_X \rightarrow X$\\
\hline
\end{tabular}
\end{center}

\vspace{.8cm}

 \begin{center}
 {\bf 3. Mod $l$ Milnor invariants of primes for $l=2,3$ }
 \end{center}
 
 In this section, we review the arithmetic of mod $2$ (resp. mod $3$) Milnor invariants of rational primes (resp. primes of $\mathbb{Q}(\zeta_3)$), which has been studied in [AMM] and [Mo1] $\sim$ [Mo4].
 \\
 \\
 {\bf 3.1. Maximal pro-$l$ Galois groups with restricted ramification for $l=2,3$.} Let $k$ be a finite algebraic number field such that $k$ contains $\zeta_l := \exp(\frac{2\pi\sqrt{-1}}{l})$  and the class number of $k$ is one.  Let $S$ be a finite subset of $s$ distinct finite primes which are not lying over $l$, $S := \{ \frak{p}_1, \dots, \frak{p}_s\}$.   Note that ${\rm N}\frak{p}_i \equiv 1$ mod $l$ ($1\leq i \leq s$). Let $k_S(l)$ denote the maximal pro-$l$  extension  of  $k$, unramified outside $S \cup S_k^{\infty}$,  in a fixed algebraic closure $\bar{k}$, where $S_k^{\infty}$ denotes the set of infinite primes of $k$. Let $G_{k,S}(l)$ denote the Galois group of $k_S(l)$ over $k$. We describe the structure of the pro-$l$ group $G_{k,S}(l)$ in a certain unobstructed case. 

 We firstly recall Iwasawa's result on the local Galois group ([Iw]). For each $i$ ($1\leq i \leq s$), let $k_{\frak{p}_i}$ be the $\frak{p}_i$-adic field with a prime element $\pi_i$. We fix an algebraic closure $\overline{k}_{\frak{p}_i}$ of $k_{\frak{p}_i}$ and an embedding $\bar{k} \hookrightarrow \overline{k}_{\frak{p}_i}$. Let $k_{\frak{p}_i}(l)$ denote the maximal pro-$l$ extension of $k_{\frak{p}_i}$ in $\overline{k}_{\frak{p}_i}$ and $G_{k_{\frak{p}_i}}(l)$ denote the Galois group of $k_{\frak{p}_i}(l)$ over $k_{\frak{p}_i}$. Then we have
$$ k_{\frak{p}_i}(l) = k_{\frak{p}_i}( \zeta_{l^a}, \sqrt[l^a]{\pi_i} \; | \; a \geq 1 ),$$
where $\zeta_{l^a}$ denotes a primitive $l^a$-th root of unity in $\bar{k}$  such that $(\zeta_{l^b})^{l^c} = \zeta_{l^{b-c}}$ for all $b \geq c$.  The local Galois group $G_{k_{\frak{p}_i}}(l)$  is  then topologically generated by the monodromy $\tau_i$ and (an extension of) the Frobenius automorphism $\sigma_i$ which are defined by
$$ \begin{array}{ll}
\tau_i(\zeta_{l^n}) := \zeta_{l^n}, & \tau_i(\sqrt[l^n]{\pi_i}) := \zeta_{l^n}\sqrt[l^n]{\pi_i},\\
\sigma_i(\zeta_{l^n}) := \zeta_{l^n}^{{\rm N}\frak{p}_i},& \sigma_i(\sqrt[l^n]{\pi_i}) := \sqrt[l^n]{\pi_i}
\end{array}
$$
and subject to the relation
$$ \tau_i^{{\rm N}\frak{p}_{i} -1}[\tau_i, \sigma_i] = 1. \leqno{(3.1.1)}$$

For each $i$ ($1\leq i \leq s$), the fixed  embedding $\bar{k} \hookrightarrow \overline{k}_{\frak{p}_i}$ gives an embedding $k_S(l) \hookrightarrow k_{\frak{p}_i}(l)$, hence a prime $\frak{P}_i$ of $k_S(l)$ lying over $\frak{p}_i$. 
We denote by the same letters $\tau_i$ and $\sigma_i$ the images of $\tau_i$ and $\sigma_i$, respectively, under the homomorphism 
$$ G_{k_{\frak{p}_i}}(l) \longrightarrow G_{k,S}(l)$$ 
induced by the embedding $k_S(l) \hookrightarrow k_{\frak{p}_i}(l)$. Then $\tau_i$ is a topological generator of the inertia group of the prime $\frak{P}_i$ and $\sigma_i$ is an extension of the Frobenius automorphism of the maximal subextension of $k_S(l)/k$ for which $\frak{P}_i$ is unramified. We call simply $\tau_i$ and $\sigma_i$  a {\it monodromy over $\frak{p}_i$} in $k_S(l)/k$ and  a {\it Frobenius automorphism over $\frak{p}_i$} in $k_S(l)/k$, respectively.

Since the ideal class group of $k$  is trivial, class field theory tells us that the monodromies $\tau_1, \dots , \tau_s$ generate topologically the global Galois group $G_{k,S}(l)$. However, they may not be a minimal set of generators  in general. In fact, noting that $k$ contains $\zeta_l$, Shafarevich's theorem ([Kc; Satz 11.8]) tells us that the minimal number $d(G_{k,S}(l))$ of generators of $G_{k,S}(l)$ is given by
 $$ d(G_{k,S}(l)) = s - r_{\mathbb{C}}(k) + \dim_{\mathbb{F}_l} B_{k,S}^{(l)}.  \leqno{(3.1.2)}$$
Here $r_{\mathbb{C}}(k)$ denotes the number of complex primes (up to conjugation) of $k$ and the obstruction $B_{k,S}^{(l)}$ is defined by 
$$B_{k,S}^{(l)} := \{ a \in k^{\times} \; | \; (a) = \frak{a}^l,  \;  a \in (k_{\frak{p}}^{\times})^l  \; \mbox{for all}\;\frak{p} \in S\cup S_k^{\infty} \}/(k^{\times})^l,  $$
where $\frak{a}$ is a fractional ideal of ${\cal O}_k$. 

In the following, we deal with the case that $l=2$ and $k = \mathbb{Q}$ or the case that $l=3$ and $k = \mathbb{Q}(\zeta_3)$. For these cases, we can determine $B_{k,S}^{(l)}$ and, moreover, we can show that the relations for minimal generators of $G_{k,S}(l)$ are given by the local relations (3.1.1).\\
\\
$\bullet$ The case that $l=2$ and $k = \mathbb{Q}$.  We have $r_{\mathbb{C}}(\mathbb{Q}) = 0$ and we can easily verify $B_{\mathbb{Q},S}^{(2)} = \{1\}$ for any $S = \{ (p_1), \dots, (p_s)\}$, where $p_i$'s are odd prime numbers. Therefore, by (3.1.2), we have $ d(G_{\mathbb{Q},S}(2)) = s$, namely, $\tau_1, \dots , \tau_s$ are minimal generators of $G_{\mathbb{Q},S}(2)$. By Koch's theorems [Kc; Satz 6.11] ([Kc; Satz 6.14]) and [Kc; Satz 11.3],  the relations for these minimal generators are given by the local relations (3.1.1). Hence, we have the following \\
\\
{\bf Theorem 3.1.3} ([Mo4; Theorem 7.4]). {\it The pro-$2$ group $G_{\mathbb{Q},S}(2)$ has the following minimal presentation
$$ \begin{array}{ll} G_{\mathbb{Q},S}(2) & = \langle \,  x_{1}, \dots , x_s  \; | \; x_{1}^{p_{1} -1}[x_{1},y_{1}] = \cdots =  x_{s}^{p_{s} -1}[x_{s},y_{s}] =1 \, \rangle \\
 & = \widehat{F_s}^{(2)}/N_S^{(2)}.
\end{array} $$
Here $\widehat{F_s}^{(2)}$ is the free pro-$2$ group generated by letters $x_1, \dots , x_s$ where $x_i$ denotes represents a monodromy $\tau_{i}$ over $p_{i}$ in $\mathbb{Q}_{S}(2)/\mathbb{Q}$, and $N_S^{(2)}$ is the closed subgroup of $\widehat{F_s}^{(2)}$ normally generated by $ x_{1}^{p_{1} -1}[x_{1},y_{1}], \dots , x_{s}^{p_{s} -1}[x_{s},y_{s}] $ where $y_{i}$ is the free pro-$2$ word in $\widehat{F_s}^{(2)}$  which represents a Frobenius automorphism over $p_{i}$ in $\mathbb{Q}_{S}(2)/\mathbb{Q}$.}\\
\\
$\bullet$ The case that $l=3$ and $k = \mathbb{Q}(\zeta_3)$. We have $r_{\mathbb{C}}(\mathbb{Q}(\zeta_3)) = 1$ and, by [AMM; Proposition 1.8], $B_{\mathbb{Q}(\zeta_3),S}^{(3)} = \{ 1\}$ if and only if $S$ contains a prime $\frak{p}$ satisfying ${\rm N}\frak{p} \equiv 4 \; \mbox{or}\; 7$ mod $9$. 
We let $S_0 := \{ \frak{p}_1, \dots , \frak{p}_{s-1}\}$ with $s \geq 2$ and ${\rm N}\frak{p}_i \equiv 1 \; \mbox{mod}\; 9$ ($1\leq i \leq s$) and let $S := S \cup \{ \frak{p}_s \}$ with ${\rm N}\frak{p}_s \equiv 4$ or $7 \; \mbox{mod}\; 9$.  By (3.1.2),  we have $ d(G_{\mathbb{Q}(\zeta_3),S}(3)) = s-1$, namely, one of $\tau_1, \dots , \tau_s$ is redundant for minimal generators of $G_{\mathbb{Q}(\zeta_3),S}(3)$. It is shown in [AMM; Proposition 1.9] that we can exclude the monodromy over $\frak{p}_s$ to obtain minimal generator of $G_{\mathbb{Q}(\zeta_3),S}(3)$.  By [Kc; Satz 6.11] ([Kc; Satz 6.14]), [Kc; Satz 11.3] and  [Kc; Satz 11.4],  we have the following\\
\\
{\bf Theorem 3.1.4} ([AMM; Theorem 1.10]). {\it  The pro-$3$ group $G_{\mathbb{Q}(\zeta_3), S}(3)$ has the following minimal presentation
$$\begin{array}{ll}  G_{\mathbb{Q}(\zeta_3), S}(3) & = \langle \,  x_{1}, \dots, x_{s-1}  \; | \; x_{1}^{{\rm N}\frak{p}_{1} -1}[x_{1},y_{1}] = \cdots = 
   x_{s}^{{\rm N}\frak{p}_{s-1} -1}[x_{s-1},y_{s-1}] =1 \, \rangle \\
 & = \widehat{F_{s-1}}^{(3)}/N_S^{(3)}.
\end{array}$$
Here   $\widehat{F_{s-1}}^{(3)}$ is the free pro-$3$ group generated by letters $x_1, \dots , x_{s-1}$ where $x_i$  represents a monodromy $\tau_{i}$ over $\frak{p}_{i}$ in $\mathbb{Q}(\zeta_3)_{S}(3)/\mathbb{Q}(\zeta_3)$, and $N_S^{(3)}$ is the closed subgroup of $\widehat{F_{s-1}}^{(3)}$ normally generated by $x_{1}^{{\rm N}\frak{p}_{1} -1}[x_{1},y_{1}], \dots,   x_{s-1}^{{\rm N}\frak{p}_{s-1} -1}[x_{s-1},y_{s-1}]$  where $y_{i}$ is the free pro-$2$ word in $\widehat{F_{s-1}}^{(3)}$  which represents a Frobenius automorphism over $\frak{p}_{i}$ in $\mathbb{Q}(\zeta_3)_{S}(3)/\mathbb{Q}(\zeta_3)$.}\\
\\
{\bf 3.2. Mod $l$ Milnor invariants of primes for $l=2,3$.} In this subsection, we recall mod $2$ Milnor invariants of rational primes and mod $3$ Milnor invariants of primes in $\mathbb{Q}(\zeta_3)$. We keep the same notations as in Subsection 3.1.\\
\\
$\bullet$ Mod $2$ Milnor invariants of rational primes.  Let $\widehat{F_s}^{(2)}$ be the free pro-$2$ group generated by $x_1, \dots , x_s$, where each $x_i$ represents a monodromy over $p_i$, as in Theorem 3.1.3. Let $\Theta_2 : \mathbb{F}_2[[\widehat{F_s}^{(2)}]]\; \stackrel{\sim}{\longrightarrow}   \mathbb{F}_2 \langle \langle X_1, \dots , X_s \rangle \rangle$ be the mod $2$ Magnus isomorphism in (1.3.2). For a multi-index $I$ and $1\leq j \leq s$, we let $\mu_2(Ij) := \mu_2(I; y_j)$, where the pro-$2$ word $y_j$ represents a Frobenius automorphism over $p_j$, so that we have 
$$ \Theta_2(y_j) = 1 + \sum_{|I| \geq 1} \mu_2(Ij) X_{I}.$$
We set $\mu_2(I) := 0$ if $|I| = 1$. Let $e_S := \max \{ e | p_1 \equiv 1 \; \mbox{mod}\; 2^e (1\leq i \leq s)\}.$\\
\\
{\bf Theorem 3.2.1} ([Mo4; 8.4]). (1) {\em For $i\neq j$, we have}
$$ \displaystyle{ (-1)^{\mu_2(ij)} = \left( \frac{p_i}{p_j} \right), }$$
{\em where $\left( \frac{p_i}{p_j} \right)$ stands for the Legendre symbol.}\\
(2) {\em Let $I = (i_1\cdots i_n)$ be a multi-index with $2 \leq n \leq 2^{e_S}$. If $\mu_2(j_1\cdots j_m) = 0$ for any proper subset $\{ j_1,\dots , j_m\}$ of $\{i_1,\dots , i_n\}$, then $\mu_2(I)$ is an invariant, called mod $2$ Milnor invariants,  of an ordered set $\{ i_1,\dots , i_n\}$.}\\
\\
$\bullet$ Mod $3$ Milnor invariants of primes in $\mathbb{Q}(\zeta_3)$. Let $\widehat{F_{s-1}}^{(3)}$ be the free pro-$3$ group generated by $x_1, \dots , x_{s-1}$, where each $x_i$ represents a monodromy over $\frak{p}_i$, as in Theorem 3.1.4. Let $\Theta_3 : \mathbb{F}_3[[\widehat{F_{s-1}}^{(3)}]]\; \stackrel{\sim}{\longrightarrow}   \mathbb{F}_3 \langle \langle X_1, \dots , X_{s-1} \rangle \rangle$ be the mod $3$ Magnus isomorphism in (1.3.2). For a multi-index $I$ and $1\leq j \leq s$, we let $\mu_3(Ij) := \mu_3(I; y_j)$, where the pro-$3$ word $y_j$ represents a Frobenius automorphism over $\frak{p}_j$, so that we have 
$$ \Theta_3(y_j) = 1 + \sum_{|I| \geq 1} \mu_3(Ij) X_{I}.$$
We set $\mu_3(I) := 0$ if $|I| = 1$. We choose the unique prime element $\pi$ of $\frak{p}_i$ ($1\leq i \leq s-1$) such that $\pi_i \equiv 1 \; \mbox{mod}\; (3\sqrt{-3})$.  \\
\\
{\bf Theorem 3.2.2.} (1) ([AMM; Theorem 3.6]). {\em For $i\neq j$, we have}
$$ \displaystyle{ \zeta_3^{\mu_3(ij)} = \left( \frac{\pi_i}{\pi_j} \right)_3,}$$
{\em where $\left( \frac{\pi_i}{\pi_j} \right)_3$ stands for the cubic residue symbol.}\\
(2) ([AMM; Proposition 4.3, Theorem 4.4]). {\em Let  $i, j, k$ be distinct indices, $1\leq i,j,k \leq s-1$. Assume that $\frak{p}_i$ and $\frak{p}_j$ are generated by rational prime numbers and that $\mu_3(ab) = 0$ for $a, b \in \{ i, j, k\}$. Then $\mu_3(ijk)$ is independent of a choice of $\frak{p}_s$ and an invariant, called the mod $3$ Milnor invariant, of an ordered set $\{ i, j, k\}$. }\\
\\
{\bf Remark 3.2.3.} As in Remark 1.3.8 (1), by using the relation between Magnus coefficients and Massey products ([Dw],[St]), it was shown in [Mo3] and [AMM; 7] that the mod $l$ Milnor invariants of primes are expressed by Massey products of the mod $l$ cohomology of the Galois group $G_{k,S}(l)$ for $l=2, k = \mathbb{Q}$ or $l=3, k= \mathbb{Q}(\zeta_3)$.
\\

\begin{center}
{\bf 4. Triple quadratic and cubic residue symbols in Ihara theory}
\end{center}

In this section, we interpret quadratic (resp. cubic) residue symbols as mod $2$ (resp. mod $3$) Milnor invariants of Galois elements in Ihara theory. \\
\\
{\bf 4.1. Triple quadratic residue symbols (R\'{e}dei symbols).}  Let $p_1$ and $p_2$ be distinct prime numbers satisfying
$$ p_i \equiv 1 \; \mbox{mod}\; 4  \; (i = 1, 2), \;\; \left( \frac{p_i}{p_j} \right) = 1 \; (1\leq i \neq j \leq 2). \leqno{(4.1.1)} $$
By the assumption (4.1.1), there are integers $x, y$ and $z$ such that 
$$x^2 - p_1y^2 -p_2z^2 = 0, (x,y,z) = 1, y \equiv 0 \; {\rm mod}\;  2, x-y \equiv 1 \; {\rm mod}\; 4. \leqno{(4.1.2)} $$
We set
$$ R^{(2)} := R^{(2)}_{\{p_1,p_2\}} := \mathbb{Q}(\sqrt{p_1}, \sqrt{p_2}, \sqrt{\alpha}),\;\;   \alpha := x + \sqrt{p_1}y.  \leqno{(4.1.3)}  $$
It is the unique Galois extension of $\mathbb{Q}$, determined by the set $\{ p_1, p_2\}$,  having the following properties: its Galois group is the dihedral group $H(\mathbb{F}_2)$ of order $8$ and it is unramified outside $p_1, p_2$ and the infinite prime ([A], [R]).  Let $p_3$ be  a prime number satisfying 
$$ p_3 \equiv 1 \; \mbox{mod}\; 4, \;\; \left( \frac{p_i}{p_j} \right) = 1 \; (1\leq i \neq j \leq 3). \leqno{(4.1.4)} $$
Let $K^{(2)} := K^{(2)}_{\{p_1,p_2\}}:= \mathbb{Q}(\sqrt{p_1}, \sqrt{p_2})$. For a prime $\tilde{\frak{p}}$  of $K^{(2)}$ lying over $p_3$, the {\em R\'{e}dei symbol} is defined by
$$ [p_1,p_2,p_3] := \frac{\left( \frac{R^{(2)}/K^{(2)}}{\tilde{\frak{p}}} \right)(\sqrt{\alpha})}{\sqrt{\alpha}}, 
 \leqno{(4.1.5)}
$$
which is independent of a choice of $\tilde{\frak{p}}$ ([A], [R]). 

By Theorems 3.2.1 applied to the case $S := \{p_1,p_2,p_3\}$ with the assumptions (4.1.1) and (4.1.4) satisfied, the mod $2$ triple Milnor invariant $\mu_2(123)$ of rational primes $\{p_1,p_2,p_3\}$ is well defined. The following theorem gives an interpretation of the R\'{e}dei symbol in terms of a mod $2$ triple Milnor invariant. \\
\\
{\bf Theorem 4.1.6} ([Mo4; 8.4]). {\em We have}
$$ \displaystyle{ (-1)^{\mu_2(123)} = [p_1,p_2,p_3]. }$$

\vspace{0.2cm}

Now we shall interpret the R\'{e}dei symbol as a mod $2$ Milnor invariant of a Galois element in Ihara theory. Following the notations in the sections 1 and 2, we consider the case where $l = 2$, $k = \mathbb{Q}$ and $A = \{ a_0, a_1, a_2, a_3 \}$ with
$$ a_0 := \infty, a_1 := 0, a_2 := x^2, a_3 := p_1 y^2.$$
Let ${\cal K}^{(2)}_{\{x, t\}}$ be the $(\mathbb{Z}/2\mathbb{Z} \times \mathbb{Z}/2\mathbb{Z})$-extension of $\mathbb{Q}(t)$ in (2.1.1) and let ${\cal R}^{(2)}_{\{x,t\}}$ be the dilogarithmic $H(\mathbb{F}_2)$-extension of $\mathbb{Q}(t)$ in (2.1.2), with $l=2$ and $c=x$:
$$  \begin{array}{ll} {\cal R}^{(2)}_{\{x,t\}} & := {\cal K}^{(2)}_{\{x, t\}}(\sqrt{\varepsilon_2(t)}) \\
                                                                      & = \mathbb{Q}(t)(\sqrt{t},\sqrt{x^2-t}, \sqrt{\varepsilon_2(t)}), \;\; \varepsilon_2(t) := x + \sqrt{t}.
\end{array} \leqno{(4.1.7)} $$
We note by (4.1.2), (4.1.3) and (4.1.7) that ${\cal K}^{(2)}_{\{x, t\}}$ and ${\cal R}^{(2)}_{\{x,t\}}$ are specialized to $K^{(2)}_{\{p_1,p_2\}}$ and $R^{(2)}_{\{p_1,p_2\}}$, respectively, by the evaluation $c = x, t = p_1y^2$:
$${\cal K}^{(2)}_{\{x,p_1y^2\}} = K^{(2)}_{\{p_1,p_2\}}, \;   {\cal R}^{(2)}_{\{x,p_1y^2\}} = R^{(2)}_{\{p_1,p_2\}}, \; \varepsilon_2(p_1y^2) = \alpha.$$
Let $\Omega_A$ be the Ihara filed of definition for $A$ in (1.1.7). By Theorem 2.6, we have
$$ R^{(2)}_{\{p_1,p_2\}} \subset \Omega_A.$$
Let ${\cal S}_A$ be as in (1.1.8). We suppose that $p_3$ satisfies $(p_3) \notin {\cal S}_A$ as well as (4.1.4). Let $\frak{P}$ be an extension of $p_3$ to $\Omega_A$ and let $\sigma_{\frak{P}}$ be the Frobenius automorphism of $\frak{P}$ over $\mathbb{Q}$.\\
\\
{\bf Proposition 4.1.8.} {\em Let the notations and assumptions be as above. For any $i, j \in \{ 1,2,3\}$, we have}
$$ \mu_2(\sigma_{\frak{P}};ij) = 0.$$
{\em Hence $\mu_2(\sigma_{\frak{P}};ij3)$ is independent of a choice of $\frak{P}$ by  Collorary 1.3.6 and so it is denoted by $\mu_2(\sigma_{p_3};ij3)$. Then we have}
$$ f_3(\sigma_{\frak{P}}) \equiv \prod_{i=1}^3 x_i^{2\mu_2(\sigma_{p_3};ii3)} \prod_{1\leq i < j \leq 3} [x_i,x_j]^{\mu_2(\sigma_{p_3};ij3)} \; \mbox{mod}\; \widehat{F_3}^{(2)}(3),$$
{\em where $\widehat{F_3}^{(2)}(3) = (\widehat{F_3}^{(2)}(2))^2[\widehat{F_3}^{(2)}, [\widehat{F_3}^{(2)},\widehat{F_3}^{(2)}]]$ is the 3rd term of the mod $2$ Zassenhaus filtration of $\widehat{F_3}^{(2)}$ (cf. (1.3.3)).}\\
\\
{\em Proof.} For  $i \neq j$, we have $a_i - a_j = \pm x^2, \pm p_1y^2, \pm p_2z^2$ by  (4.1.2). So the first assertion follows from Theorem 1.3.7 and the  assumptions (4.1.1) and (4.1.4), and so we have $ f_3(\sigma_{\frak{P}}) \in  \widehat{F_3}^{(2)}(2).$ Note by (1.3.3) that $\widehat{F_3}^{(2)}(2)/\widehat{F_3}^{(2)}(3)$ has  a basis $x_1^2, x_2^2, x_3^2, [x_1,x_2], [x_2, x_3], [x_1,x_3]$ over $\mathbb{F}_2$. Then the second assertion follows from the definition of mod $2$ Milnor invariants. $\;\; \Box$
\\
\\
{\bf Proposition 4.1.9.}  {\em The monodromy transformation of $\sqrt{\varepsilon_2(t)}$ along the pro-$2$ longitude $f_3(\sigma_{\overline{\frak{p}}})$  is given by}
$$ \sqrt{\varepsilon_2(t)}\;  \mapsto  \; [p_1,p_2,p_3] \sqrt{\varepsilon_2(t)}.$$ 
\\
{\em Proof.} First, we note the followings.\\
(i) By induction on $n \geq 0$, we easily see
$$ \frac{d^n}{dt^n} \sqrt{\varepsilon_2(t)} \,\Big|_{t= p_1y^2} = \sqrt{\alpha} \cdot \Phi_n(\alpha, y\sqrt{p_1})$$
for some $\Phi_n(X_1,X_2) \in \mathbb{Q}(X_1,X_2)$.\\
(ii) By (4.1.4) and (4.1.5), we have 
$$\chi_2(\sigma_{\frak{P}}) \equiv 1 \,\mbox{mod}\; 4, \; \sigma_{\frak{P}}(\sqrt{p_1}) = \sqrt{p_1}, \; \mbox{and}\; \sigma_{\frak{P}}(\sqrt{\alpha}) = [p_1,p_2,p_3] \sqrt{\alpha}.$$
(iii) We easily see $\iota_0(\sqrt{\varepsilon_2(t)}) \in \mathbb{Q}\{ \{ 1/t \} \}$.

Then the monodromy transformation of $\sqrt{\varepsilon_2(t)}$ along the pro-$2$ longitude $f_3(\sigma_{\frak{P}}) =  s_0(\sigma_{\frak{P}}) \cdot \gamma_3^{-1}\cdot s_3(\sigma_{\frak{P}})^{-1} \cdot \gamma_3$ 
is given as follows:
$$ \begin{array}{ll}
\iota_0(\sqrt{\varepsilon_2(t)}) & \stackrel{\gamma_3}{\longrightarrow} \displaystyle{ \sqrt{\alpha} \sum_{n=0}^{\infty} \Phi_n(\alpha, y\sqrt{p_1}) (t- p_1y^2)^n} \;\; (\mbox{by (i)})\\
                                    & \stackrel{s_3(\sigma_{\frak{P}})^{-1}}{\longrightarrow} \displaystyle{ [p_1,p_2,p_3] \sqrt{\alpha} \sum_{n=0}^{\infty} \Phi_n(\alpha, y\sqrt{p_1}) (t- p_1y^2)^n} \; \; (\mbox{by (ii)}) \\
                                     & \stackrel{\gamma_3^{-1}}{\longrightarrow}  [p_1,p_2,p_3]  \iota_0(\sqrt{\varepsilon_2(t)}) \\
                                     & \stackrel{s_0(\sigma_{\frak{P}})}{\longrightarrow} [p_1,p_2,p_3]  \iota_0(\sqrt{\varepsilon_2(t)})  \; \; (\mbox{by (iii)}). \;\; \Box\\
                                     \end{array}
                                     $$ 
\\
{\bf Theorem 4.1.10.} {\em We have}
$$ [p_1,p_2,p_3]  = (-1)^{\mu_2(\sigma_{p_3};123)}$$
{\em and hence} 
$$ \mu_2(123) = \mu_2(\sigma_{p_3};123).$$
\\
{\em Proof.} Since ${\cal K}^{(2)}_{\{x, t\}}$ is a metabelian extension of $k(t)$, any element of $\widehat{F_3}^{(2)}(3)$ acts on $\sqrt{\varepsilon_2(t)}$ trivially. Then the assertion follows from Proposition 4.1.8, Corollary 2.1.6 and Proposition 4.1.9. $\;\; \Box$\\
\\
{\bf 4.2. Triple cubic residue symbols.} Let $\frak{p}_1$ and $\frak{p}_2$ be distinct primes of $\mathbb{Q}(\zeta_3)$ with ${\rm N}\frak{p}_i \equiv 1$ mod $9$. We assume that each $\frak{p}_i$ is generated by a rational prime number. Let $\pi_i$ be the unique prime element of $\frak{p}_i$ such that $\pi_i \equiv \; 1 \; \mbox{mod}\; (3\sqrt{-3})$.  We assume that 
$$ \left( \frac{\pi_i}{\pi_j} \right)_3 = 1 \;\; (1\leq i \neq j \leq 2). \leqno{(4.2.1)} $$ 
Let $K_1 := k(\sqrt[3]{\pi_1})$ and let $\tau$ be the generator of the Galois group of $k(\sqrt[3]{\pi_1})$ defined by $\tau(\sqrt[3]{\pi_1}) = \zeta_3 \sqrt[3]{\pi_1}$. By (4.2.1), there is $\alpha$ in ${\cal O}_{K_1}$ such that 
$${\rm N}_{K_1/k}(\alpha) = \pi_2 z^3, (\alpha) = \frak{P}^e \frak{B}^f \; {\rm with}\;  (e,3) = 1,\;  (\frak{B}, 3) = 1, \;  f \equiv 0 \; \mbox{mod} \; 3, \leqno{(4.2.2)}$$
where $z \in \mathbb{Z}[\zeta_3]$  and $\frak{P}, \frak{B}$ are ideals of $\mathbb{Z}[\zeta_3]$. We let
$$\theta := \tau(\alpha)(\tau^2(\alpha))^2 \leqno{(4.2.3)}$$
and set
$$ R^{(3)} := R^{(3)}_{\{\frak{p}_1, \frak{p}_2\}} := \mathbb{Q}(\zeta_3)(\sqrt[3]{\pi_1}, \sqrt[3]{\pi_2}, \sqrt[3]{\theta}). \leqno{(4.2.4)} $$
It is the unique Galois extension of $\mathbb{Q}(\zeta_3)$, determined by the set $\{ \frak{p}_1, \frak{p}_2\}$,  having the following properties: its Galois group is isomorphic to $H(\mathbb{F}_3)$ and only $\frak{p}_1$ and $\frak{p}_2$ are ramified with ramification indices being $3$ ([AMM; Theorem 5.11, Corollary 5.12]). Let $\frak{p}_3$ be  a prime of $\mathbb{Q}(\zeta_3)$ such that ${\rm N}\frak{p}_3 \equiv 1 \; \mbox{mod}\; 9$ and let $\pi_3$ be the unique prime element in $\frak{p}_3$ such that $\pi_3 \equiv \; 1 \; \mbox{mod}\; (3\sqrt{-3})$. We assume that 
$$ \left( \frac{\pi_i}{\pi_j} \right)_3 = 1 \; (1\leq i \neq j \leq 3). \leqno{(4.2.5)} $$
Let $K^{(3)} := K^{(3)}_{\{\frak{p}_1, \frak{p}_2\}}:= \mathbb{Q}(\zeta_3)(\sqrt[3]{\pi_1}, \sqrt[3]{\pi_2})$. For a prime $\tilde{\frak{p}}$  of $K^{(3)}$ lying over $\frak{p}_3$, we define the {\em triple cubic residue symbol} by
$$ [\frak{p}_1, \frak{p}_2, \frak{p}_3]_3 := \frac{\left( \frac{R^{(3)}/K^{(3)}}{\tilde{\frak{p}}} \right)(\sqrt[3]{\theta})}{\sqrt[3]{\theta}}, 
 \leqno{(4.2.6)}
$$
which is independent of the choice of a prime $\tilde{\frak{p}}$. 

By Theorems 3.2.2 applied to the case $S := \{\frak{p}_1, \frak{p}_2, \frak{p}_3\}$ with the assumptions (4.2.1) and (4.2.5) satisfied, the mod $3$ triple Milnor invariant $\mu_3(123)$ of  primes $\{\frak{p}_1, \frak{p}_2, \frak{p}_3\}$ is well defined. The following theorem gives an interpretation of the R\'{e}dei symbol in terms of a mod $3$ triple Milnor invariant. \\
\\
{\bf Theorem 4.2.7} ([AMM; Definition 6.2, Theorem 6.3]). {\em We have}
$$ \displaystyle{ \zeta_3^{\mu_3(123)} = [\frak{p}_1, \frak{p}_2, \frak{p}_3]_3. }$$

\vspace{.08cm}

Now we shall interpret the triple cubic residue symbol as a mod $3$ Milnor invariant of a Galois element in Ihara theory. In the following, we assume that $\alpha \in {\cal O}_{K_1}$ in (4.2.2)  is of the form
$$ \alpha = x + y \sqrt[3]{\pi_1} \leqno{(4.2.8)} $$
for some $x, y \in \mathbb{Q}(\zeta_3)$ and so ${\rm N}_{K_1/k}(\alpha) = \pi_2z^2$ in (4.2.2) and $\theta$ in (4.2.3) are written as
$$ x^3 + \pi_1y^3 = \pi_2 z^3. \leqno{(4.2.9)}$$
and $$ \theta = (x+\zeta_3 y\sqrt[3]{\pi_1})(x+\zeta_3^2 y\sqrt[3]{\pi_1})^2.  \leqno{(4.2.10)}$$
 Following the notations in the sections 1 and 2, we consider the case where $l = 3$, $k = \mathbb{Q}(\zeta_3)$ and $A = \{ a_0, a_1, a_2, a_3 \}$ with
$$ a_0 := \infty, a_1 := 0, a_2 := x^3, a_3 := -\pi_1 y^3.$$
Let ${\cal K}^{(3)}_{\{x, t\}}$ be the $(\mathbb{Z}/3\mathbb{Z} \times \mathbb{Z}/3\mathbb{Z})$-extension of $\mathbb{Q}(\zeta_3)(t)$ in (2.1.1) and let ${\cal R}^{(3)}_{\{x,t\}}$ be the dilogarithmic $H(\mathbb{F}_3)$-extension of $\mathbb{Q}(\zeta_3)(t)$ in (2.1.2), with $l=3$ and $c=x$:
$$  \begin{array}{ll} {\cal R}^{(3)}_{\{x,t\}} & := {\cal K}^{(3)}_{\{x, t\}}(\sqrt[3]{\varepsilon_3(t)}) \\
                                                                       & = \mathbb{Q}(\zeta_3)(t)(\sqrt[3]{t},\sqrt[3]{x^3-t}, \sqrt[3]{\varepsilon_3(t)}), \;\; \varepsilon_3(t) := (x - \zeta_3 \sqrt[3]{t})(x - \zeta_3^2 \sqrt[3]{t})^2. \end{array} \leqno{(4.2.11)}$$
We note by (4.2.4), (4.2.9), (4.2.10) and (4.2.11) that ${\cal K}^{(3)}_{\{x, t\}}$ and ${\cal R}^{(3)}_{\{x,t\}}$ are specialized to $K^{(3)}_{\{\frak{p}_1, \frak{p}_2\}}$ and $R^{(3)}_{\{\frak{p}_1,\frak{p}_2\}}$, respectively, by the evaluation $c = x, t = -\pi_1y^3$:
$${\cal K}^{(3)}_{\{x, -\pi_1y^3\}} = K^{(3)}_{\{\frak{p}_1, \frak{p}_2\}}, \;   {\cal R}^{(3)}_{\{x, -\pi_1y^3\}} = R^{(3)}_{\{\frak{p}_1,\frak{p}_2\}}, \; \varepsilon_3(-\pi_1y^3) =\theta.$$
Let $\Omega_A$ be the Ihara filed of definition in (1.1.7). By Theorem 2.6, we have
$$ R^{(3)}_{\{\frak{p}_1, \frak{p}_2\}} \subset \Omega_A.$$
Let ${\cal S}_A$ be as in (1.1.8). We suppose that $\frak{p}_3$ satisfies $\frak{p}_3 \notin {\cal S}_A$ as well as (4.2.5). Let $\frak{P}$ be an extension of  $\frak{p}_3$ to $\Omega_A$ and let $\sigma_{\frak{P}}$ be the Frobenius element of $\frak{P}$ over $\mathbb{Q}(\zeta_3)$.\\
\\
{\bf Proposition 4.2.12.} {\em Let the notations and assumptions be as above. For any $i, j \in \{ 1,2,3\}$, we have}
$$ \mu_3(\sigma_{\frak{P}};ij) = 0.$$
{\em Hence $\mu_3(\sigma_{\frak{P}};ij3)$ is independent of a choice of $\frak{P}$ by  Corollary 1.3.6 and so it is denoted by $\mu_3(\sigma_{\frak{p}_3};ij3)$. Then we have}
$$ f_3(\sigma_{\frak{P}}) \equiv \prod_{1\leq i < j \leq 3} [x_i,x_j]^{\mu_3(\sigma_{\frak{p}_3};ij3)} \; \mbox{mod}\; \widehat{F_3}^{(3)}(3),$$
{\em where $\widehat{F_3}^{(3)}(3) = (\widehat{F_3}^{(3)})^3[\widehat{F_3}^{(3)}, [\widehat{F_3}^{(3)},\widehat{F_3}^{(3)}]]$ is the 3rd term of the mod $3$ Zassenhaus filtration of $\widehat{F_3}^{(3)}$ (cf. (1.3.3)).}\\
\\
{\em Proof.} For $i \neq j$, we have $a_i - a_j = \pm x^3, \pm \pi_1y^3, \pm \pi_2 z^3$ by (4.2.9). So the first assertion follows from Theorem 1.3.7 and the  assumptions (4.2.1) and (4.2.5), and so $ f_3(\sigma_{\frak{P}}) \in \widehat{F_3}^{(3)}(2).$ Note by (1.3.3) that $\widehat{F_3}^{(3)}(2)/\widehat{F_3}^{(3)}(3)$ has a basis $[x_1, x_2], [x_2,x_3], [x_1,x_3]$ over $\mathbb{F}_3$. Then the second assertion follows from the definition of mod $3$ Milnor invariants.  $\;\; \Box$
\\
\\
{\bf Proposition 4.2.13.} {\em The monodromy transformation of $\sqrt[3]{\varepsilon_3(t)}$ along the pro-$3$ longitude $f_3(\sigma_{\frak{P}})$  is given by}
$$ \sqrt[3]{\varepsilon_3(t)}\;  \mapsto  \; [\frak{p}_1, \frak{p}_2, \frak{p}_3]^{-1} \sqrt[3]{\varepsilon_3(t)}.$$ 
\\
{\em Proof.}  The proof goes in a way similar to that of Proposition 4.1.9. First, we note the followings.\\
(i) By induction on $n \geq 0$, we easily see
$$ \frac{d^n}{dt^n} \sqrt[3]{\varepsilon_3(t)} \Big|_{t= -\pi_1y^3} = \sqrt[3]{\theta} \cdot \Phi_n(\theta, y\sqrt[3]{\pi_1})$$
for some $\Phi_n(X_1,X_2) \in \mathbb{Q}(X_1,X_2)$.\\
(ii) By (4.2.5) and (4.2.6), we have 
$$  \chi_3(\sigma_{\frak{P}}) \equiv 1 \; \mbox{mod}\; 9, \;  \sigma_{\frak{P}}(\sqrt[3]{\pi_1}) = \sqrt[3]{\pi_1}, \; \mbox{and}\; \sigma_{\frak{P}}(\sqrt[3]{\theta}) = [\frak{p}_1, \frak{p}_2,\frak{p}_3]_3 \sqrt[3]{\theta}.$$
(iii) We easily see $\iota_0(\sqrt[3]{\varepsilon_3(t)}) \in \mathbb{Q}(\zeta_3)\{ \{ 1/t \} \}$.

Then the monodromy transformation of $\sqrt[3]{\varepsilon_3(t)}$ along the pro-$3$ longitude $f_3(\sigma_{\frak{P}}) =  s_0(\sigma_{\frak{P}}) \cdot \gamma_3^{-1}\cdot s_3(\sigma_{\frak{P}})^{-1} \cdot \gamma_3$ 
is given as follows:
$$ \begin{array}{ll}
\iota_0(\sqrt[3]{\varepsilon_3(t)}) & \stackrel{\gamma_3}{\longrightarrow} \displaystyle{ \sqrt[3]{\theta} \sum_{n=0}^{\infty} \Phi_n(\theta, y\sqrt[3]{\pi_1}) (t+ \pi_1y^3)^n} \; \; (\mbox{by (i)})\\
                                    & \stackrel{s_3(\sigma_{\frak{P}})^{-1}}{\longrightarrow} \displaystyle{ [\frak{p}_1,\frak{p}_2,\frak{p}_3]_3^{-1} \sqrt[3]{\theta} \sum_{n=0}^{\infty} \Phi_n(\theta, y\sqrt[3]{\pi_1}) (t+ \pi_1y^3)^n} \;\; (\mbox{by (ii)}) \\
                                     & \stackrel{\gamma_3^{-1}}{\longrightarrow}  [\frak{p}_1, \frak{p}_2, \frak{p}_3]_3^{-1}  \iota_0(\sqrt[3]{\varepsilon_3(t)}) \\
                                     & \stackrel{s_0(\sigma_{\frak{P}})}{\longrightarrow}   [\frak{p}_1, \frak{p}_2, \frak{p}_3]_3^{-1}   \iota_0(\sqrt[3]{\varepsilon_3(t)}) \;\; (\mbox{by (iii)}). \;\; \Box\\
                                     \end{array}
                                     $$ 
\\
{\bf Theorem 4.2.14.} {\em We have}
$$ [\frak{p}_1, \frak{p}_2, \frak{p}_3]_3^{-1}  = \zeta_3^{\mu_3(\sigma_{\frak{p}_3};123)}$$
{\em and hence} 
$$ - \mu_3(123) = \mu_3(\sigma_{\frak{p}_3};123).$$
\\
{\em Proof.} Since ${\cal K}^{(3)}_{\{x, t\}}$ is a metabelian extension of $k(t)$, any element of $\widehat{F_3}^{(3)}(3)$ acts on $\sqrt[3]{\varepsilon_3(t)}$ trivially. Then the assertion follows from Proposition 4.2.12, Corollary 2.1.6 and Proposition 4.2.13. $\;\; \Box$\\
\\
{\bf Example 4.2.15.} The assumptions (4.2.1) and (4.2.8) are satisfied for the cases $(-\pi_1, -\pi_2) = (17,53), (17, 467), (107, 449), (431, 233)$ etc.(This computation is due to Y. Mizusawa.)
\\
Let $(-\pi_1, -\pi_2) = (17,53)$. Then we can take $x = 8, y = 3, z = -1$  and so 
$$\alpha = 8 - 3\sqrt[3]{17}, \; \theta = (8 - 3\zeta_3\sqrt[3]{17})(8 - 3\zeta_3^2\sqrt[3]{17})^2$$
and 
$$ {\cal R}^{(3)}_{\{8, 3^3\cdot 17\}} = R^{(3)}_{\{(17),(53)\}} = \mathbb{Q}(\zeta_3)(\sqrt[3]{17},\sqrt[3]{53}, \sqrt[3]{\theta}). $$
By [AMM; Example 6.4], for $-\pi_3 = 71, 89, 107, 179, 197$, we have
$$ \begin{array}{l} \mu_3(\sigma_{(71)};123) = - \mu_3(123) = 1,  \; \mu_3(\sigma_{(89)};123) = - \mu_3(123) = 2,\\
 \mu_3(\sigma_{(107)};123) = - \mu_3(123) = 1,\;  \mu_3(\sigma_{(179)};123) = - \mu_3(123) = 2,  \\
\mu_3(\sigma_{(197)};123) = - \mu_3(123) = 2.  \end{array}$$ 
\\
{\small 
\begin{flushleft}
{\bf References}\\
{[A]} F. Amano, On R\'{e}dei's dihedral extension and triple reciprocity law, Proc. Japan Acad. Ser. A Math. Sci.  {\bf 90}  (2014),  no. 1, 1--5.\\
{[AMM]} F. Amano, Y. Mizusawa, M. Morishita, On mod $3$  triple Milnor invariants and triple cubic residue symbols in the Eisenstein number field, Research in Number Theory  {\bf 4 }, no. 1, Art. 7, 29 pp, 2018. \\
{[AI]} G. Anderson, Y. Ihara, Pro-$l$ branched coverings of $\mathbb{P}^1$ and higher circular $l$-units, Ann. of Math. (2)  {\bf 128}  (1988),  no. 2, 271--293.\\
{[DDMS]} J.D. Dixon, M.P.F. du Sautoy, A. Mann, D. Segal, Analytic pro-$p$ Groups, second edition, Cambridge Stud. Adv. Math., vol.{\bf 61}, Cambridge University Press, Cambridge, 1999.\\
{[Be]} A. Beilinson, Higher regulators and values of $L$-functions, J. Soviet Math., {\bf 30} (1985), 2036--2070.\\
{[Bl]} S. Bloch, The dilogarithm and extensions of Lie algebras,  Algebraic K -theory, Evanston 1980 (Proc. Conf., Northwestern Univ., Evanston, Ill., 1980),  pp. 1--23, Lecture Notes in Math., {\bf 854}, Springer, Berlin-New York, 1981.\\
{[Br1]} J.-L. Brylinski, Loop spaces, characteristic classes and geometric quantization, Progress in Mathematics, {\bf 107}. Birkhauser Boston, Inc., Boston, MA, 1993. \\
{[Br2]} J.-L. Brylinski, Holomorphic gerbes and the Beilinson regulator, $K$-theory (Strasbourg, 1992). 
Ast\'{e}risque  No. {\bf 226}  (1994), 8, 145-174. \\
{[Dw]} W. G.  Dwyer, Homology, Massey products and maps between groups, J. Pure Appl. Algebra  {\bf 6}  (1975), no. 2, 177--190.\\
{[De]} P. Deligne, Le symbole mod\'{e}r\'{e}, Inst. Hautes Etudes Sci. Publ. Math.  No. {\bf 73}  (1991), 147--181.\\
{[Gi]} J. Giraud, Cohomologie non ab\'{e}lienne, Die Grundlehren der Mathematischen Wissenschaften, {\bf 179}. Springer-Verlag, New York-Berlin, 1971. \\
{[Gr]} A. Grothendieck,  Rev\^{e}tements \'{E}tales et Groupe Fondamental, Lecture Notes in Mathematics, {\bf 224}, Springer-Verlag, 1971.\\
{[I1]} Y. Ihara, Profinite braid groups, Galois representations and complex multiplications, Ann. of Math. (2)  {\bf 123}  (1986),  no. 1, 43--106.\\
{[I2]} Y. Ihara, Arithmetic analogues of braid groups and Galois representations,  Braids (Santa Cruz, CA, 1986),  245--257, Contemp. Math., {\bf 78}, Amer. Math. Soc., Providence, RI, 1988.\\
{[I3]} Y. Ihara, Braids, Galois groups, and some arithmetic functions,  Proceedings of the International Congress of Mathematicians, Vol. I, II (Kyoto, 1990),  99--120, Math. Soc. Japan, Tokyo, 1991.\\
{[Iw]} K. Iwasawa, On Galois groups of local fields, Trans. Amer. Math. Soc.  {\bf 80}  (1955), 448--469. \\
{[Kc]} H. Koch, Galoissche Theorie der $p$-Erweiterungen, Springer-Verlag, Berlin-New York; VEB Deutscher Verlag der Wissenschaften, Berlin, 1970.\\
{[Kd]} H. Kodani, Arithmetic topology on braid and absolute Galois groups, Doctoral Thesis, Kyushu University, 2017.\\
{[KMT]} H. Kodani, M. Morishita, Y. Terashima, Arithmetic topology in Ihara theory,  Publ. Res. Inst. Math. Sci., Kyoto Univ., {\bf 53}, No.4. (2017), 629--688.\\
{ [Mi1] } J. Milnor, Link groups, Ann. of Math. {\bf 59} (1954), 177--195.\\
{[Mi2]} J. Milnor, Isotopy of links, in Algebraic Geometry and Topology, A symposium in honor of S. Lefschetz ( edited by R.H. Fox, D.C. Spencer and A.W. Tucker), 280--306 Princeton University Press, Princeton, N.J., 1957.\\
{ [Mo1] } M. Morishita, Milnor's link invariants attached to certain Galois groups over {\bf Q}, Proc. Japan Acad. Ser, A {\bf 76} (2000), 18--21.\\
{ [Mo2] } M. Morishita, On certain analogies between knots and primes, J. Reine Angew. Math. {\bf 550} (2002), 141--167.\\
{ [Mo3] } M. Morishita, Milnor invariants and Massey products for prime numbers, Compos. Math., {\bf 140} (2004), 69--83. \\
{[Mo4]} M. Morishita, Knots and Primes -- An Introduction to Arithmetic Topology, Universitext, Springer, 2011.\\
{[MK]} K. Murasugi, B.I. Kurpita, A study of braids, Mathematics and its Applications, {\bf 484}, Kluwer Academic Publishers, Dordrecht, 1999. \\
{[N]} H. Nakamura, Tangential base points and Eisenstein power series,  Aspects of Galois theory (Gainesville, FL, 1996),  202--217, 
London Math. Soc. Lecture Note Ser., {\bf 256}, Cambridge Univ. Press, Cambridge, 1999. \\
{[NW]} H. Nakamura, Z. Wojtkowiak,  On explicit formulae for $l$-adic polylogarithms,  In: Arithmetic fundamental groups and noncommutative algebra (Berkeley, CA, 1999),  285--294, Proc. Sympos. Pure Math., 70, Amer. Math. Soc., Providence, RI, 2002.\\
{ [R] } L. R\'{e}dei, Ein neues zahlentheoretisches Symbol mit Anwendungen auf die Theorie der quadratischen Zahlk\"{o}rper I, J. Reine Angew. Math., {\bf 180} (1939), 1-43.\\
{[So]} C. Soul\'{e}, $K$-th\'{e}orie des anneaux d'entiers de corps de nombres et cohomologie \'{e}tale,  Invent. Math.  {\bf 55}  (1979), no. 3, 251--295.\\
{[St] } D. Stein, Massey products in the cohomology of groups with applications to link theory, Trans.
Amer. Math. Soc. {\bf 318}  (1990), 301--325.\\
{ [T] } V. Turaev, The Milnor invariants and Massey products, (Russian) Studies in topology, II. Zap. Nau\v{c}n. Sem. Leningrad. Otdel. Mat. Inst. Steklov. (LOMI)  {\bf 66},  (1976), 189--203, 209--210. \\
{[W1]} Z. Wojtkowiak, On $l$-adic iterated integrals, I. Analog of Zagier conjecture, Nagoya Math. J.  {\bf 176}  (2004), 113--158.\\
{[W2]} Z. Wojtkowiak, On $l$-adic iterated integrals, II. Functional equations and $l$-adic polylogarithms, Nagoya Math. J.  {\bf 177}  (2005), 117--153.\\
{[W3]} Z. Wojtkowiak, On $l$-adic iterated integrals, III. Galois actions on fundamental groups. Nagoya Math. J.  {\bf 178}  (2005), 1--36.\\
{[W4]} Z. Wojtkowiak, On $l$-adic iterated integrals, IV. Ramification and generators of Galois actions on fundamental groups and torsors of paths,  Math. J. Okayama Univ. {\bf  51}  (2009), 47--69.\\
{[W5]} Z. Wojtkowiak, A remark on nilpotent polylogarithmic extensions of the field of rational functions of one variable over $\mathbb{C}$, Tokyo J. Math.  {\bf 30}  (2007),  no. 2, 373--382.\\ 
\end{flushleft} 
\vspace{.5cm}
Hikaru Hirano\\
Faculty of Mathematics, Kyushu University \\
744, Motooka, Nishi-ku, Fukuoka, 819-0395, JAPAN \\
e-mail: ma218019@math.kyushu-u.ac.jp\\
\\
Masanori Morishita\\
Faculty of Mathematics, Kyushu University \\
744, Motooka, Nishi-ku, Fukuoka, 819-0395, JAPAN \\
e-mail: morisita@math.kyushu-u.ac.jp
}
\end{document}